\documentclass[12pt]{amsart}

\usepackage{amssymb,amsmath}
\usepackage{amsthm}
\usepackage{tikz}
\usepackage{color}
\usepackage{fullpage}
\usepackage{booktabs}
\usepackage{enumerate}
\usepackage[all]{xy}

\usetikzlibrary{patterns}
\usetikzlibrary{calc}

\newcommand{\eps}{\epsilon}

\newcommand{\bbA}{\mathbb{A}}
\newcommand{\bbC}{\mathbb{C}}
\newcommand{\bbD}{\mathbb{D}}

\newcommand{\bbM}{\mathbb{M}}

\newcommand{\bbR}{\mathbb{R}}

\newcommand{\bbS}{\mathbb{S}}
\newcommand{\bbZ}{\mathbb{Z}}

\newcommand{\calD}{\mathcal{D}}
\newcommand{\calE}{\mathcal{E}}
\newcommand{\calF}{\mathcal{F}}
\newcommand{\calG}{\mathcal{G}}

\newcommand{\calM}{\mathcal{M}}
\newcommand{\calN}{\mathcal{N}}

\newcommand{\calP}{\mathcal{P}}
\newcommand{\calQ}{\mathcal{Q}}

\newcommand{\calS}{\mathcal{S}}

\newcommand{\calW}{\mathcal{W}}

\newcommand{\frakD}{\mathfrak{D}}
\newcommand{\frakL}{\mathfrak{L}}
\newcommand{\frakR}{\mathfrak{R}}

\newcommand{\fraka}{\mathfrak{a}}
\newcommand{\frakh}{\mathfrak{h}}

\newcommand{\End}{\mathrm{End}}

\newcommand{\Hom}{\mathrm{Hom}}

\newcommand{\diag}{\mathrm{diag}}

\newcommand{\ord}{\mathrm{ord}}
\newcommand{\sign}{\mathrm{sign}}

\newcommand{\lag}{\mathfrak{g}}
\newcommand{\lah}{\mathfrak{h}}

\newcommand{\lak}{\mathfrak{k}}

\newcommand{\lam}{\mathfrak{m}}

\newcommand{\las}{\mathfrak{s}}
\newcommand{\lat}{\mathfrak{t}}

\newcommand{\laso}{\mathfrak{so}}

\newcommand{\laZ}{\mathfrak{Z}}

\newcommand{\Spin}{\mathrm{Spin}}

\newcommand{\sph}{\mathrm{sph}}

\newcommand{\Rad}{\mathrm{Rad}}

\newcommand{\ns}{\mathrm{ns}}
\newcommand{\ct}{\mathrm{ct}}

\newcommand{\rFs}[5]{\,_{#1}F_{#2} \left( \genfrac{.}{.}{0pt}{}{#3}{#4}
\ ;#5 \right)}

\renewenvironment{proof}[1][]{\noindent{\scshape Proof{#1}. }}{\qed}
\newenvironment{proofempty}[1][]{\noindent{\scshape Proof{#1}. }}{}

\theoremstyle{plain}
\newtheorem{theorem}{Theorem}[section]
\newtheorem{lemma}[theorem]{Lemma}
\newtheorem{proposition}[theorem]{Proposition}
\newtheorem{corollary}[theorem]{Corollary}

\newtheorem{definition}[theorem]{Definition}

\theoremstyle{definition}
\newtheorem{example}[theorem]{Example}
\newtheorem{remark}[theorem]{Remark}
\newtheorem{notation}[theorem]{Notation}

\title{Non-symmetric Jacobi polynomials of type $BC_1$ as vector-valued polynomials Part 2:\\ Shift operators}

\date{\today}

\author[van Horssen]{M. van Horssen}
\author[van Pruijssen]{M. van Pruijssen}

\address[van Horssen]{KU Leuven, Department of Mathematics, Celestijnenlaan 200B, 3001 Leuven, Belgium}
\email{max.vanhorssen@kuleuven.be}
\address[van Pruijssen]{Radboud University, IMAPP-Mathematics, Heyendaalseweg 135, 6525 AJ NIJMEGEN, the Netherlands}
\email{m.vanpruijssen@math.ru.nl}


\begin{document}
\begin{abstract}
We study non-symmetric Jacobi polynomials of type $BC_1$ by means of vector-valued and matrix-valued orthogonal polynomials. The interpretation as matrix-valued orthogonal polynomials allows us to introduce shift operators for the non-symmetric Jacobi polynomials. The shift operators are differential-reflection operators and we present four of these operators that are fundamental in the sense that they generate all shift operators. Moreover, the symmetrizations of these fundamental shift operators are the fundamental shift operators for the symmetric Jacobi polynomials of type $BC_1$.

For the realization of non-symmetric Jacobi polynomials of type $BC_1$ as invariant $\bbC^2$-valued Laurent polynomials, we introduce a homomorphism that is analogous to the Harish-Chandra homomorphism for the symmetric Jacobi polynomials of type $BC_1$. For geometric root multiplicities, the non-symmetric Jacobi polynomials of type $BC_1$ can be interpreted as spherical functions and we show that our Harish-Chandra homomorphism in this context is related to the Lepowsky homomorphism via the radial part map.
\end{abstract}

\maketitle


\section{Introduction}\label{S:Intro}
This paper is the continuation of the study the non-symmetric Jacobi polynomials of type $BC_1$ by means of vector-valued and matrix-valued orthogonal polynomials \cite{vHvP1}, and its main purpose is to introduce and classify shift operators for these polynomials in all their forms. We also study the related classification problem of determining all differential operators that have the non-symmetric Jacobi polynomials of type $BC_1$, viewed as invariant vector-valued Laurent polynomials, as eigenfunctions.

The symmetric Jacobi polynomials for general root systems have been introduced by Heckman and Opdam in \cite{HeOp87}, where they are defined using the Gram–Schmidt process with respect to a partial ordering. A priori, two polynomials with incomparable labels need not be orthogonal. The symmetric Jacobi polynomials are simultaneous eigenfunctions of a family of commuting differential operators, which can be used to prove their orthogonality. Originally, the main tools to study these operators were the Heckman operator, the trigonometric version of the Dunkl operator, and the Harish-Chandra homomorphism. The surjectivity of the Harish-Chandra homomorphism was conjectured by Heckman and Opdam in \cite{HeOp87}, and it was proven by Opdam in \cite{Opd88b} by means of analytic methods. An important question within Heckman-Opdam theory is to determine and relate the several normalizations, i.e., the leading coefficient, the norms, and the evaluation at the identity. The determination of the norms of the symmetric Jacobi polynomials is closely related to the well-known Macdonald constant term conjectures \cite{Mac82}. For the infinite families of root systems these conjectures have been settled by Gunson, Wilson, and Macdonald, see \cite{Gun62, Wil62, Mac82}. In \cite{Opd89}, Opdam gave a uniform proof using shift operators for the symmetric Jacobi polynomials, and thereby proving the constant term conjectures for the remaining exceptional cases. The existence of these shift operators had been established by Opdam in \cite{Opd88b} by analytic methods, and their algebraic properties were studied before in \cite{Opd88a}. Later Heckman showed the existence by more elementary (algebraic) means in \cite{Hec91}.

The Heckman operators are equivariant for the action of the Weyl group of the root system, however, they fail to commute with each other. In \cite{Che91}, Cherednik introduced an alternative trigonometric version of the Dunkl operator, known as the Cherednik operator. These operators commute with one another but they are no longer Weyl group equivariant. In this way, Cherednik established a pivotal connection between symmetric Jacobi polynomials and the (graded) Hecke algebra. The commutation of these operators opened the door to study their simultaneous eigenfunctions, which gave rise to the non-symmetric Jacobi polynomials \cite{Opd95}. The symmetric Jacobi polynomials can be obtained by symmetrizing the non-symmetric Jacobi polynomials, so the latter can be regarded as more fundamental objects. For a more complete historical account we refer the reader to \cite{Opd00}. The surjectivity of the Harish-Chandra can also be obtained from this Hecke algebra interpretation, see \cite[Thm.2.10]{Opd95}. The norms of the non-symmetric Jacobi polynomials have been calculated in \cite{Opd95} by studying certain representations of the Hecke algebra, and by reducing them to the norms of the symmetric Jacobi polynomials. In \cite{Opd95}, the norms of the symmetric Jacobi polynomials are calculated using the so-called shift principle \cite[Lem.5.2]{Opd95}. This principle can also be used to prove the existence of the shift operators using the Cherednik operator, see \cite{Opd00}.

A natural question is whether there exist shift operators for the non-symmetric Jacobi polynomials, and in that case if these are more fundamental in the sense that we recover the shift operators for the symmetric Jacobi polynomials after some symmetrization action.

We establish the existence of shift operators for the non-symmetric Jacobi polynomials of type $BC_1$ in Section \ref{S:ShiftNonSym}. In the symmetric case, the cornerstone of the existence of shift operators, in the algebraic construction, is the shift principle. There is no analog of this principle available to us in the non-symmetric case. However, for the non-symmetric Jacobi polynomials of type $BC_1$ we circumvent this issue by using a decomposition of these polynomials in terms of their symmetric counterparts. This was established in \cite{vHvP1} by realizing the non-symmetric and symmetric Jacobi polynomials as matrix-valued orthogonal polynomials which are  related by a simple transformation that relies on the fact that the matrix-weights are conjugated. We recall the construction of this result in Section \ref{S:prerequisites}, and we utilize this in Section \ref{S:shiftMVOP} to construct shift operators for the matrix-valued orthogonal polynomials. More precisely, we construct four matrix-valued shift operators that are fundamental in the sense that an arbitrary matrix-valued shift operator is a product of these fundamental shift operators and some differential operator that has the matrix-valued orthogonal polynomials as eigenfunctions. This is one of the main results of this paper that is called the structure theorem for the matrix-valued shift operators. Furthermore, we also classify all differential operators that have the matrix-valued orthogonal polynomials as eigenfunctions in Section \ref{subsection:DW}. The matrix-valued differential operator that corresponds to the Cherednik operator plays a crucial role in this classification. Previously, this classification has been found by Castro and Grünbaum in \cite{CaGr06} for the Chebyshev case, where the scalar part of the matrix-weight is given by a Chebyshev weight.

The matrix-valued shift operators are transferred to differential-reflection operators that are shift operators for the non-symmetric Jacobi polynomials. These operators are characterized by satisfying a transmutation property with the Cherednik operator, and we call such operators non-symmetric shift operators. Similarly, we will refer to shift operators for the symmetric Jacobi polynomials as symmetric shift operators. There is a correspondence between non-symmetric and matrix-valued shift operators, and the structure theorem for the non-symmetric shift operators is a direct consequence of this. In Section \ref{S:ShiftNonSym}, we develop some theory for the non-symmetric shift operators that is in analogy with the theory for symmetric shift operators (of type $BC_1$) as discussed in \cite[\S3.3]{HeSc94}. This allows us to prove the structure theorem for the non-symmetric shift operators without appealing to its matrix-valued version, which relies on the structure theorem for the symmetric shift operators. As an application of this theory, we demonstrate in Section \ref{S:normalization} how to use the non-symmetric shift operators to calculate the normalizations of the non-symmetric Jacobi polynomials, from which we recover the known formulas of type $BC_1$ from \cite{Opd95} more directly.

The non-symmetric Jacobi polynomials can also be interpreted as $\bbC^2$-valued Laurent polynomials that are invariant for an action of the Weyl group, see \cite[\S2]{vHvP1}. In Section \ref{S:HC}, we study the algebra of differential operators that act on the space of $\bbC^2$-valued Laurent polynomials which admit a so-called (matrix-valued) constant term. The matrix-valued differential operator that encodes the KZ-equations \cite[\S3]{Opd95}, which also corresponds to the Cherednik operator, lies in this algebra and plays a pivotal role. We show that the matrix-valued constant term map is injective on the commutant of this differential operator. By introducing a $\rho$-shift in the matrix-valued constant term map, we obtain a matrix-valued analog of the Harish-Chandra homomorphism that was considered by Heckman and Opdam in \cite{HeOp87} for the symmetric Jacobi polynomials. If we furthermore impose some Weyl group invariance on the differential operators, then the matrix-valued Harish-Chandra homomorphism restricts to an isomorphism onto a polynomial algebra in one variable, cf. \cite[Thm.1.2.9]{HeSc94}. Furthermore, the $\bbC^2$-valued version of the non-symmetric Jacobi polynomials are eigenfunctions for the elements of this subalgebra, and the matrix-valued Harish-Chandra homomorphism describes the eigenvalues.

For geometric root multiplicities, the $\bbC^2$-valued version of the non-symmetric Jacobi polynomials of type $BC_1$ have been identified with spherical functions for the compact symmetric pair $(\Spin(2m+2), \Spin(2m+1))$ in \cite{vHvP1}. In this case, the algebra of differential operators that have the spherical functions as eigenfunctions is the image of the radial part map, which is defined on a commutative subquotient of the universal enveloping algebra $U(\laso(2m+2, \bbC))$. The eigenvalues are related to the Lepowsky homomorphism, and this ultimately allows us to connects the matrix-valued Harish-Chandra homomorphism with the Lepowsky homomorphism in Section \ref{S:SF}.

Finally, we note that the matrix-valued version of the non-symmetric Jacobi polynomials of type $BC_1$ with root multiplicity $(0,\nu)$ can be identified with the matrix-valued Gegenbauer polynomials of size $2\times 2$ that have been studied in \cite{MR3735699}. Here differentiation with respect to the variable is identified as a shift operator, and it coincides with one of our fundamental matrix-valued shift operators. While preparing this manuscript, shift operators for the non-symmetric Askey-Wilson polynomials have been found in \cite{vHS}, which we can regard as the $q$-analog of our non-symmetric shift operators because they are recovered for a suitable $q\to1$ limit.


\section{Prerequisites}\label{S:prerequisites}
We recall some facts and notation from \cite[\S2]{vHvP1}. The root system $R = \{ \pm\eps, \pm2\eps \} \subset \bbR^*$ of type $BC_1$ has Weyl group $\bbZ_2 = <s>$ and weight lattice $P=\bbZ\eps$. The group algebra of $P$ is identified with the algebra of Laurent polynomials $\bbC[z^{\pm 1}]$ and admits an action of the Weyl group, given by $s(f)(z)=f(z^{-1})$. The root multiplicity is a $\bbZ_2$-invariant function $k : R \to \bbC$ that is represented by $k = (k_1, k_2)$, where $k_1$ is the value on $\eps$ and $k_2$ on $2\eps$. The space $\bbC[z^{\pm 1}]$ is equipped with a sesquilinear pairing
\begin{equation}\label{eq:pairing_on_s1}
(f, g)_k = \int_{S^1} \overline{f(z)}g(z)\delta_k(z)\frac{dz}{iz},
\end{equation}
where $\delta_k(z) = \left(1 - (z + z^{-1}) / 2 \right)^{k_1}\left(1 - (z^2 + z^{-2}) / 2\right)^{k_2}$. For $k_1, k_2 \geq 0$ the pairing $(\cdot, \cdot)_k$ is an inner product. Applying the Gram-Schmidt process on the ordered basis $(1, z, z^{-1}, z^2, \ldots)$ yields an orthogonal basis $(E(n,k) \mid n \in \bbZ)$ of $\bbC[z^{\pm 1}]$ such that
\begin{itemize}
\item $E(n, k) = z^n +$ lower order terms,
\item $(E(n, k), z^m)_k = 0$ for all monomials $z^m < z^n$.
\end{itemize}
The Laurent polynomials $E(n,k)$ are called the non-symmetric Jacobi polynomials of type $BC_{1}$. 
The differential-reflection operator
\begin{equation}\label{eq:cher_op}
D_k = z\partial_z + k_1\frac{1}{1 - z^{-1}}(1 - s) + 2k_2\frac{1}{1 - z^{-2}}(1 - s) - \rho(k),
\end{equation}
where $\rho(k) = \tfrac{1}{2}(k_1 + 2k_2)$, is called the Cherednik operator. The operator $D_k$ acts on $\bbC[z^{\pm 1}]$, it is formally self-adjoint for the pairing $(\cdot, \cdot)_k$, and it has the polynomials $E(n, k)$ as eigenfunctions
\begin{equation}\label{eq:BC1_eigenvalue}
\begin{aligned}
D_k E(n,k) &= (n + \rho(k))E(n,k), & n > 0, \\
D_k E(n,k) &= (n - \rho(k))E(n,k), & n \leq 0.
\end{aligned}
\end{equation}
Let $x = (z + z^{-1}) / 2$ and recall that $\bbC[x] = \bbC[z^{\pm 1}]^{\bbZ_2} \subset \bbC[z^{\pm 1}]$ is the ring of invariants. The pairing $(\cdot, \cdot)_k$ restricted to $\bbC[z^{\pm 1}]^{\bbZ_2}$ and restricted to $\bbC[x]$ is given by
\begin{equation*}
(f, g)_k = \int_{-1}^1\overline{f(x)}g(x)w_k(x)dx,
\end{equation*}
where $w_k(x) = 2^{k_1 + 2k_2}(1 - x)^{k_1 + k_2 - \frac{1}{2}}(1 + x)^{k_2 - \frac{1}{2}}$. Application of the Gram-Schmidt process on the ordered basis $(1, x, x^2, \ldots)$ with respect to the pairing $(\cdot, \cdot)_k$ yields an orthogonal basis $(p(n, k) \mid n \geq 0)$ of $\bbC[x]$ such that
\begin{itemize}
\item $p(n, k) = (z + z^{-1})^n +$ lower order terms,
\item $(p(n, k), (z+z^{-1})^m)_k = 0$ for all $m < n$.
\end{itemize}
The polynomials $p(n,k)$ are called symmetric Jacobi polynomials of type $BC_{1}$. 
With the convention $\alpha = k_1 + k_2 - \tfrac{1}{2}$, $\beta = k_2 - \tfrac{1}{2}$, we have
\begin{equation}\label{eq:monic_p_to_clas_p}
p(n, k) = \frac{2^{2n} n!}{(n + \alpha + \beta + 1)_n}P_n^{(\alpha, \beta)}((z + z^{-1}) / 2),
\end{equation}
where
\begin{equation*}
P^{(\alpha, \beta)}_n(x) = \frac{(\alpha + 1)_n}{n!}\rFs{2}{1}{-n, n + \alpha + \beta + 1}{\alpha + 1}{\frac{1 - x}{2}}
\end{equation*}
are the classical Jacobi polynomials.


\subsection{Shift operators for the symmetric Jacobi polynomials}\label{S:shiftops_for_sym_jac}
Shift operators for the polynomials $p(n, k)$ allow us to relate the polynomials $p(n, k)$ with polynomials belonging to a family with a different parameter $k'$. With these operators the normalizations of the polynomials $p(n, k)$ can be determined. We will refer to these operators as symmetric shift operators to differentiate them from the shift operators for the polynomials $E(n, k)$, i.e., the non-symmetric shift operators. We recall from \cite{HeSc94, Opd88a} some well-known facts about symmetric shift operators for the root system of type $BC_1$.

Let $T$ be a differential-reflection operator that restricts to an operator on $\bbC[x]$. Then there exists a differential operator $\beta(T)$ on $\bbC[z^{\pm 1}]$ that coincides with $T$ on $\bbC[x]$, see \cite[Thm.2.12(ii)]{Opd95}. The operator $D_k^2$ commutes with the action of $\bbZ_2$ on $\bbC[z^{\pm 1}]$, so it restricts to an operator on $\bbC[x]$. Hence we can define $ML_k = \beta(D_k^2)$, which is called the modified Laplacian, and given explicitly by
\begin{equation}\label{eq:mod_lap}
ML_k = (z\partial_z)^2 + k_1 \frac{1 + z^{-1}}{1 - z^{-1}} z\partial_z + 2 k_2 \frac{1 + z^{-1}}{1 - z^{-1}} z\partial_z + (\tfrac{1}{2}k_1 + k_2)^2.
\end{equation}
The operator $ML_k$ has the polynomials $p(n, k)$ as eigenfunctions,
\begin{equation*}
ML_k p(n, k) = (n + \tfrac{1}{2}k_1 + k_2)^2 p(n, k), \quad n \geq 0.
\end{equation*}

Let $\Delta = (z - z^{-1})$ and denote the localization of $\bbC[z^{\pm 1}]$ along $\Delta$ by $\bbC_\Delta[z^\pm]$.

\begin{definition}\label{def:sym_shift_op}
An element $S(k) \in \bbC_\Delta[z^\pm] \otimes \bbC[z\partial_{z}]$ is called a symmetric shift operator with shift $\ell \in 2\bbZ \times \bbZ$ if
\begin{enumerate}[(i)]
    \item $S(k)$ satisfies the transmutation property, i.e., $S(k) ML_k = ML_{k + \ell} S(k)$,
    \item $S(k)$ acts as a differential operator on $\bbC[x]$.
\end{enumerate}
\end{definition}

\begin{remark}\label{rem:sym_shiftops_weyl_algebra}
It should be noted that Definition \ref{def:sym_shift_op} differs from \cite[Def.2.1]{Opd88a} and \cite[Def.3.1.1]{HeSc94}. We take the coefficients as in \cite[Def.3.1.1]{HeSc94}, but for fixed $k$, and instead of requiring an asymptotic expansion, either as a formal or convergent power series, we assume Definition \ref{def:sym_shift_op}(ii), which is implied by the former in conjunction with the transmutation property, see \cite[Prop.2.2]{Opd88a}. Finally, a symmetric shift operator, or any element of $\bbC_\Delta[z^\pm] \otimes \bbC[z\partial_{z}]$ satisfying Definition \ref{def:sym_shift_op}(ii), can be viewed as an element of the Weyl algebra $\bbC[x, \partial_x]$, see \cite[Prop.2.2]{Opd88a}.
\end{remark}

\begin{remark}
The assumption in Definition \ref{def:sym_shift_op} that $\ell \in 2\bbZ \times \bbZ$ is not a restriction, because for $\ell \notin 2\bbZ \times \bbZ$ there are no nonzero differential operators that satisfy the transmutation property, see \cite[Cor.3.9(b)]{Opd88a}.
\end{remark}

It is well-known that the operator $\partial_x$, called the forward shift, is a shift operator for the Jacobi polynomials $P^{(\alpha, \beta)}_n(x)$, see \cite[(9.8.7)]{KLS10}, which is then also a shift operator for the polynomials $p(n, k)$ by \eqref{eq:monic_p_to_clas_p}. This operator can be extended to an operator on $\bbC[z^{\pm 1}]$ that is a symmetric shift operator in the sense of Definition \ref{def:sym_shift_op}. 

\begin{definition}[{\cite[Prop.3.3.1]{HeSc94}}]
The forward and backward shift operators are given by
\begin{equation*}
G_+(k) = \frac{1}{z - z^{-1}} z\partial_z, \quad G_-(k) = (z - z^{-1}) z\partial_z + (k_1 + 2k_2 - 1)(z + z^{-1}) + 2k_1, 
\end{equation*}
which have shifts $(0, 1)$ and $(0, -1)$, respectively. The contiguity shift operators are
\begin{equation*}
E_+(k) = \frac{1 + z^{-1}}{1 - z^{-1}} z\partial_z + k_2 - \tfrac{1}{2}, \quad E_-(k) = \frac{1 - z^{-1}}{1 + z^{-1}} z\partial_z + k_1 + k_2 - \tfrac{1}{2},
\end{equation*}
which have shifts $(2, -1)$ and $(-2, 1)$, respectively. These four operators are called the fundamental symmetric shift operators.
\end{definition}

\begin{lemma}[{\cite[Cor.3.1.4, Cor.3.3.2]{HeSc94}}]\label{lem:action_fund_sym_shiftsops}
Any fundamental symmetric shift operator $S(k)$ with shift $\ell$ acts on the polynomials $p(n, k)$ by
\begin{equation*}
S(k)p(n, k) = \eta(n, S(k)) p(n - \rho(\ell), k + \ell), \quad n \geq 0,
\end{equation*}
with the convention $p(-1, k) = 0$. The so-called shift factors are given by
\begin{align*}
\eta(n, G_+(k)) &= n, \qquad \qquad \; \; \; \eta(n, G_-(k)) = n + k_1 + 2k_2 - 1, \\
\eta(n, E_+(k)) &= n + k_2 - \tfrac{1}{2}, \quad \eta(n, E_-(k)) = n + k_1 + k_2 - \tfrac{1}{2}.
\end{align*}
\end{lemma}

For each $\ell \in 2\bbZ \times \bbZ$, we denote the space of symmetric shift operators with shift $\ell$ by $\bbS(\ell, k)$. Note that the composition of $S(k) \in \bbS(\ell, k)$ and $S'(k + \ell) \in \bbS(\ell', k + \ell)$ gives $S'(k + \ell) S(k)$ $\in \bbS(\ell + \ell', k)$. It follows that $\bbS(\ell, k)$ is a right $\bbS(0, k)$-module, and that $\bbS(0, k)$ is an algebra. By composing the fundamental symmetric shift operators, we can see that the spaces $\bbS(\ell, k)$ are non-empty.

\begin{proposition}[{\cite[Cor.3.3.3]{HeSc94}}]\label{prop:general_g}
Let $\ell = \varepsilon_1 N_1 (2, -1) + \varepsilon_2 N_2 (0, 1) \in 2\bbZ \times \bbZ$ with $N_1, N_2 \in \bbZ_{\geq 0}$ and $\varepsilon_1, \varepsilon_2 \in \{ \pm 1 \}$. Define the operator
\begin{equation*}
G(\ell, k) = \prod_{i = 0}^{N_1 - 1} G(\varepsilon_1 (2, -1), k + \varepsilon_1 i(2, -1) + \varepsilon_2 N_2(0, 1)) \prod_{j = 0}^{N_2 - 1} G(\varepsilon_2 (0, 1), k + \varepsilon_2 j(0, 1)),
\end{equation*}
with the initial conventions
\begin{align*}
G((0, 1), k) &= G_+(k), \quad G((0, -1), k) = G_-(k), \\
G((2, -1), k) &= E_+(k), \quad\thinspace G((-2, 1), k) = E_-(k).
\end{align*}
Then $G(\ell, k)$ is a symmetric shift operator with shift $\ell$, i.e., $G(\ell, k) \in \bbS(\ell, k)$.
\end{proposition}

For later reference we record the following property of the operators $G(\ell, k)$.

\begin{lemma}[{\cite[(3.4)]{Opd89}}]\label{lem:sym_g_adjoint}
For $\ell \in 2\bbZ \times \bbZ$, we have
\begin{equation*}
(G(\ell, k)f, g)_{k + \ell} = (f, G(-\ell, k + \ell)g)_k, \quad f,g \in \bbC[z^{\pm 1}].
\end{equation*}
\end{lemma}

Finally, the structure of the spaces $\bbS(\ell, k)$ has been determined by Opdam in \cite{Opd88a}. 

\begin{theorem}[Structure theorem for $\bbS(\ell, k)$, {\cite[Cor.3.9]{Opd88a}}]\label{thm:sym_shift_str_thm}\text{ }
\begin{enumerate}[(i)]
    \item The algebra $\bbS(0, k)$ is generated by the modified Laplacian $ML_k$.
    \item The $\bbS(0, k)$-module $\bbS(\ell, k)$ is of rank one and freely generated by $G(\ell, k)$.
    \item The generators satisfy the relations
    \begin{equation*}
    G(\ell + \ell', k) = G(\ell', k + \ell) G(\ell, k) = G(\ell, k + \ell') G(\ell', k)
    \end{equation*}
    whenever $\sign(\ell_1) = \sign(\ell'_1)$ and $\sign(\frac{1}{2}\ell_1 + \ell_2) = \sign(\frac{1}{2}\ell'_1 + \ell'_2)$. 
\end{enumerate}
\end{theorem}


\subsection{$\bbC^2$-valued Jacobi polynomials}
We recall the definition of the vector-valued Jacobi polynomials and their orthogonality from \cite[\S2.2]{vHvP1}. For later reference we record the result.

\begin{lemma}\label{lem:steinberg_basis}
Writing $x = (z + z^{-1})/2$, we have $\bbC[z^{\pm 1}] = \bbC[x] \oplus \bbC[x]z$. Moreover, if $f = f_1 + f_2 z$, with $f_1, f_2 \in \bbC[x]$, then
\begin{equation*}
f_1 = \frac{z s(f) - z^{-1} f}{z - z^{-1}}, \quad f_2 = \frac{f - s(f)}{z - z^{-1}}.
\end{equation*}
\end{lemma}

\begin{proof}
For the first statement see \cite[\S2.2]{vHvP1}. The second statement follows directly from the observation $s(f) = f_1 + f_2 z^{-1}$.
\end{proof}

The space $\bbC[x] \otimes \bbC^2$ is a $\bbC[x]$-module by component-wise multiplication, and it is equipped with the sesquilinear pairing
\begin{equation}\label{eq:vv_pairing}
(\calP, \calQ)_k = \int_{-1}^1 \calP(x)^* \calW_k(x) \calQ(x) dx, \quad \calW_k(x) = \begin{pmatrix} 1 & x \\ x & 1 \end{pmatrix} w_k(x),
\end{equation}
that is an inner product for $k_1, k_2 \geq 0$. The $\bbC[x]$-module isomorphism
\begin{equation*}
\Upsilon : \bbC[z^{\pm 1}] \to \bbC[x] \otimes \bbC^2, \; f \mapsto \underline{f} = (f_1, f_2)^T
\end{equation*}
respects the sesquilinear forms. Here $(f_1, f_2)^T$ denotes a column vector.

Define $\calP(n, k) = \underline{E(n, k)}$ for $n \in \bbZ$, which agrees with \cite[Def.2.3]{vHvP1}. These polynomials are called the vector-valued Jacobi polynomials. We conclude that $(\calP(n, k) \mid n \in \bbZ)$ is an orthogonal basis of $\bbC[x] \otimes \bbC^2$.


\subsection{$\bbM_2$-valued Jacobi polynomials}
There are two closely related versions of the matrix-valued Jacobi polynomials, see \cite[\S2.3]{vHvP1}. We state their definition and give the relation between them.

Let $\bbM_2 = \bbC^{2\times2}$ be the algebra of $2\times 2$-matrices with entries in $\bbC$, and denote the algebra of $\bbM_2$-valued polynomials by $\bbM_2[x]$. The $\bbM_2$-valued sesquilinear pairing for the space $\bbM_2[x]$ is given by using the formula \eqref{eq:vv_pairing} for two polynomials $\calP, \calQ \in \bbM_2[x]$. For $k_1, k_2 \geq 0$ this defines an $\bbM_2$-valued inner product. Define $\calM(N, k) = (\calP(-N, k) \; \calP(N + 1, k))$ for $N \geq 0$, which is the first version of matrix-valued Jacobi polynomials. Note that $(\calM(N, k) \mid N \geq 0)$ is an orthogonal basis of $\bbM_2[x]$ with respect to $(\cdot, \cdot)_k$.

Recall from \cite[\S2.1]{vHvP1} that $\calW_k$ diagonalizes,
\begin{equation*}
U\calW_kU^* = \diag(w_{k_+}, w_{k_-}), \quad U = \frac{1}{\sqrt{2}}\begin{pmatrix} 1 & -1 \\ 1 & 1 \end{pmatrix},
\end{equation*}
with $k_+ = k + (1, 0)$ and $k_- = k + (-1, 1)$. Define another $\bbM_2$-valued sesquilinear pairing
\begin{equation}\label{eq:mv_inner_prod_prime}
(\calP, \calQ)'_k = \int_{-1}^1 \calP(x)^* \diag(w_{k_+}, w_{k_-}) \calQ(x)dx = (\calP_1, \calQ_1)_{k_+} + (\calP_2, \calQ_2)_{k_-},
\end{equation}
for which $(\calN(N, k) \mid N \geq 0)$, with $\calN(N,k) = \diag(p(N, k_+), p(N, k_-))$, is a natural orthogonal basis of $\bbM_2[x]$.

\begin{proposition}\label{prop:E_P_rel}
The two families of matrix-valued Jacobi polynomials $\{\calN(N,k) \}_{N \geq 0}$ and $\{ \calM(N,k) \}_{N \geq 0}$ are related by
\begin{equation*}
\calM(N,k) = U^{-1} \calN(N,k)UC(N,k), \quad C(N,k) = \begin{pmatrix} 1 & c_N(k) \\ 0 & 1 \end{pmatrix},
\end{equation*}
where $c_N(k) = \frac{k_1}{1 + 2N + k_1 + 2k_2}$. Here $C(N, k)$ is the leading term of  $\calM(N,k)$, which is invertible.
\end{proposition}

\begin{proofempty}
See \cite[Prop.2.5]{vHvP1} and adjust for the new normalization of $\calN(N,k)$, i.e.,
\pushQED{\qed}
\begin{equation*}
\calN(N,(\alpha,\beta)) = \diag(P^{(\alpha + 1, \beta)}_N, P^{(\alpha, \beta + 1)}_N) = \frac{(N + \alpha + \beta + 2)_N}{2^{2N}N!}\calN(N, k). \qedhere
\end{equation*}
\popQED
\end{proofempty}

Following \cite[\S2.3]{vHvP1}, the polynomials $\calM(N, k)$ are eigenfunctions of the operator
\begin{equation*}
\calD_k =
\begin{pmatrix}
-x & -1 \\
1 & x
\end{pmatrix}
\partial_x +
\begin{pmatrix}
-\rho(k) & k_1 \\
0 & 1 + \rho(k)
\end{pmatrix},
\end{equation*}
in the following sense,
\begin{equation}\label{eq:BC1DIffOpM}
\calD_k \calM(N, k) = \calM(N,k)\Lambda(N, k), \quad \Lambda(N, k) =
\diag(-N - \rho(k),  N + 1 + \rho(k)).
\end{equation}
Here $\partial_x$ is component-wise differentiation with respect to $x$. It follows that the operator
\begin{equation}\label{BC1:frakD}
\frakD_k = U\calD_k U^{-1} =
\begin{pmatrix}
\frac{1}{2}(1 - k_1) & -E_+(k_-) \\
-E_-(k_+) & \frac{1}{2}(1 + k_1)
\end{pmatrix}
\end{equation}
has the polynomials $\calN(N, k)$ as eigenfunction, i.e.,
\begin{equation*}
\frakD_k\calN(N,k)=\calN(N,k)\frakL(N,k), \quad \frakL(N,k) =
\begin{pmatrix}
\frac{1}{2}(1 - k_1) & -N - k_2 - \frac{3}{2} \\
-N - k_1 - k_2 - \frac{3}{2} & \frac{1}{2}(1 + k_1)
\end{pmatrix}.
\end{equation*}


\section{$\bbM_2$-valued shift operators}\label{S:shiftMVOP}


\subsection{Fundamental shift operators for the $\bbM_2$-valued Jacobi polynomials}
We construct four matrix-valued shift operators that act on the polynomials $\calN(N, k)$, and using Proposition \ref{prop:E_P_rel} we transfer these to matrix-valued shift operators acting on the polynomials $\calM(N, k)$. Moreover, we show that the corresponding shift factors are diagonal, which allows us to restrict the action to the columns of $\calM(N, k)$.

\begin{definition}\label{def:mv_shift_op}
An element $\widehat{S}(k) \in \bbC_\Delta[z^\pm] \otimes \bbC[z\partial_z] \otimes \bbM_2$ is called a (matrix-valued) $\frakD_k$-shift operator with shift $\ell$ if
\begin{enumerate}[(i)]
    \item $\widehat{S}(k)$ satisfies the transmutation property, i.e., $\widehat{S}(k) \frakD_k = \frakD_{k + \ell} \widehat{S}(k)$,
    \item $\widehat{S}(k)$ acts as a matrix-valued differential operator on $\bbM_2[x]$.
\end{enumerate}
We define (matrix-valued) $\calD_k$-shift operators similarly.
\end{definition}

\begin{remark}\label{rem:mv_shiftops_weyl_algebra}
We can view a matrix-valued shift operator, or any element of $\bbC_\Delta[z^\pm] \otimes \bbC[z\partial_{z}] \otimes \bbM_2$ satisfying Definition \ref{def:mv_shift_op}(ii), as an element of the Weyl algebra $\bbC[x, \partial_x] \otimes \bbM_2$, which follows from applying Remark \ref{rem:sym_shiftops_weyl_algebra} entry-wise.
\end{remark}

\begin{proposition}\label{prop:frakD_shift_ops}
Let $\gamma_+, \gamma_- \in \bbC$, and denote the differential operators
\begin{equation*}
\begin{pmatrix} G_\pm(k_+) & 0 \\ 0 & G_\pm(k_-) \end{pmatrix}, \quad \begin{pmatrix} E_+(k_+) & 0 \\ \gamma_+ & E_+(k_-) \end{pmatrix}, \quad 
\begin{pmatrix} E_-(k_+) & \gamma_- \\ 0 & E_-(k_-) \end{pmatrix}
\end{equation*}
by $\widehat{G}_\pm(k)$, $\widehat{E}_+(k; \gamma_+)$, and $\widehat{E}_-(k; \gamma_-)$, respectively. Then $\widehat{G}_+(k)$ and $\widehat{G}_-(k)$ are $\frakD_k$-shift operators, and $\widehat{E}_+(k; \gamma_+)$ and $\widehat{E}_-(k; \gamma_-)$ are $\frakD_k$-shift operators if and only if $\gamma_\pm = -1$. The shifts of $\widehat{G}_\pm(k)$ and $\widehat{E}_\pm(k) = \widehat{E}_\pm(k; -1)$ are given by $(0, \pm 1)$ and $(\pm 2, \mp 1)$, respectively.
\end{proposition}

\begin{proof}
Using the relations from Theorem \ref{thm:sym_shift_str_thm}(iii), we find that
\begin{align*}
\widehat{G}_\pm(k) \frakD_k &= \begin{pmatrix} \frac{1}{2}(1 - k_1)G_\pm(k_+) & -G_\pm(k_+)E_+(k_-) \\ -G_\pm(k_-)E_-(k_+) & \frac{1}{2}(1 + k_1)G_\pm(k_-) \end{pmatrix} \\
&= \begin{pmatrix} \frac{1}{2}(1 - k_1)G_\pm(k_+) & -E_+(k_- + (0, \pm 1))G_\pm(k_-) \\ -E_-(k_+ + (0, \pm 1))G_\pm(k_+) & \frac{1}{2}(1 + k_1)G_\pm(k_-) \end{pmatrix} = \frakD_{k + (0, \pm 1)} \widehat{G}_+(k).
\end{align*}
We compute
\begin{equation*}
\widehat{E}_+(k; \gamma_+) \frakD_k = \begin{pmatrix} \frac{1}{2}(1 - k_1)E_+(k_+) & -E_+(k_+)E_+(k_-) \\ -E_+(k_-)E_-(k_+) + \frac{1}{2}(1 - k_1)\gamma_+ & (\frac{1}{2}(1 + k_1) - \gamma_+)E_+(k_-) \end{pmatrix}
\end{equation*}
and
\begin{equation*}
\frakD_{k + (2, -1)} \widehat{E}_+(k; \gamma_+) = \begin{pmatrix} (\frac{1}{2}(-1 - k_1) - \gamma_+)E_+(k_+) & -E_+(k_+)E_+(k_-) \\ -E_-(k_+ + (2, -1))E_+(k_+) + \frac{1}{2}(k_1 + 3)\gamma_+ & \frac{1}{2}(3 + k_1)E_+(k_-) \end{pmatrix}.
\end{equation*}
Using the relation
\begin{equation}\label{eq:EplusEmin_rel}
E_+(k + (-2, 1))E_-(k) = E_-(k + (2, -1))E_+(k) + k_1,
\end{equation}
which can be verified by a small computation, we see that
$\widehat{E}_+(k; \gamma_+) \frakD_k = \frakD_{k + (2, -1)} \widehat{E}_+(k; \gamma_+)$ if and only if $\gamma_+ = -1$. Similarly, we compute
\begin{align*}
\widehat{E}_-(k; \gamma_-) \frakD_k &= \begin{pmatrix} (\frac{1}{2}(1 - k_1) - \gamma_-)E_-(k_+) & -E_-(k_+)E_+(k_-) + \frac{1}{2}(1 + k_1) \gamma_- \\ -E_-(k_-)E_-(k_+) & \frac{1}{2}(1 + k_1) E_-(k_-) \end{pmatrix}, \\
\frakD_{k + (-2, 1)} \widehat{E}_-(k; \gamma_-) &= \begin{pmatrix} \frac{1}{2}(3 - k_1)E_-(k_+) & -E_+(k_- + (-2, 1))E_-(k_-) + \frac{1}{2}(3 - k_1) \gamma_- \\ -E_-(k_-)E_-(k_+) & (\frac{1}{2}(-1 + k_1) - \gamma_-) E_-(k_-) \end{pmatrix},
\end{align*}
and see that  $\widehat{E}_+(k; \gamma_-) \frakD_k = \frakD_{k + (-2, 1)}\widehat{E}_+(k; \gamma_-)$ if and only if $\gamma_- = -1$ by using \eqref{eq:EplusEmin_rel}.
\end{proof}

\begin{corollary}\label{cor:calD_shift_ops}
The differential operators $\widehat{G}_\pm^U(k) = U^{-1}\widehat{G}_\pm(k)U$, $\widehat{E}_\pm^U(k) = U^{-1}\widehat{G}_\pm(k)U$ are $\calD_k$-shift operators.
\end{corollary}

We call the operators from Proposition \ref{prop:frakD_shift_ops} (or Corollary \ref{cor:calD_shift_ops}) the fundamental $\calD_k$- (or $\frakD_k$-) shift operators. For stating the action of these fundamental matrix-valued shift operators, it will be convenient to introduce the conventions $\calN(-1, k) = 0$ and $\calM(-1, k)= 0$.

\begin{lemma}\label{lem:matrix_G_shift_P}
If $\widehat{S}(k)$ is a fundamental $\frakD_k$-shift operator with shift $\ell$, then
\begin{equation*}
\widehat{S}(k) \calN(N,k) = \calN(N - \rho(\ell), k + \ell) H(N, \widehat{S}(k)), \quad N \geq 0,
\end{equation*}
with shift factors
\begin{align*}
H(N, \widehat{G}_+(k)) &= N, \qquad H(N, \widehat{G}_-(k)) = N + k_1 + 2k_2, \\
H(N, \widehat{E}_+(k)) &= \begin{pmatrix} N + k_2 - \tfrac{1}{2} & 0 \\ -1 & N + k_2 + \tfrac{1}{2} \end{pmatrix}, \\
H(N, \widehat{E}_-(k)) &= \begin{pmatrix} N + k_1 + k_2 + \tfrac{1}{2} & -1 \\ 0 & N + k_1 + k_2 - \tfrac{1}{2} \end{pmatrix}.
\end{align*}
\end{lemma}

\begin{proof}
This follows from Lemma \ref{lem:action_fund_sym_shiftsops}.
\end{proof}

\begin{proposition}\label{prop:matrix_G_shift_E}
If $\widehat{S}^U(k)$ is a fundamental $\calD_k$-shift operator with shift $\ell$, then
\begin{equation*}
\widehat{S}^U(k) \calM(N,k) = \calM(N - \rho(\ell), k + \ell) H(N, \widehat{S}^U(k)), \quad N \geq 0,
\end{equation*}
with (diagonal) shift factors
\begin{align*}
H(N, \widehat{G}^U_+(k)) &= N, \qquad H(N, \widehat{G}^U_-(k)) = N + k_1 + 2k_2, \\
H(N, \widehat{E}^U_+(k)) &= \diag(N + k_2 - \tfrac{1}{2}, N + k_2 + \tfrac{1}{2}), \\
H(N, \widehat{E}^U_-(k)) &= \diag(N + k_1 + k_2 - \tfrac{1}{2}, N + k_1 + k_2 + \tfrac{1}{2}).
\end{align*}
\end{proposition}

\begin{proof}
Using Proposition \ref{prop:E_P_rel} and Lemma \ref{lem:matrix_G_shift_P}, we find that
\begin{equation*}
\widehat{G}^U_+(k) \calM(N, k) = \calM(N - 1, k + (0, 1)) C(N - 1, k + (0, 1))^{-1} U^{-1} H(N, \widehat{G}_+(k)) UC(N, k)
\end{equation*}
for all $N \geq 1$, and
\begin{multline*}
\widehat{G}^U_-(k) \calM(N, k) =\\
\calM(N + 1, k + (0, -1)) C(N + 1, k + (0, -1))^{-1} U^{-1} H(N, \widehat{G}_-(k)) UC(N, k)
\end{multline*}
for all $N \geq 0$. Note that $C(N \mp 1, k + (0, \pm 1)) = C(N, k)$ and that $H(N, \widehat{G}_-(k))$ is a scalar multiple of the identity. This implies that
\begin{align*}
H(N, \widehat{G}^U_+(k)) &= H(N, \widehat{G}_+(k)), \quad N \geq 1, \\
H(N, \widehat{G}^U_-(k)) &= H(N, \widehat{G}_-(k)), \quad N \geq 0,
\end{align*}
which yield the shift factors. The case $H(0, \widehat{G}^U_+(k)) = 0$ can be verified directly.

Similarly, we find that
\begin{equation*}
\widehat{E}^U_\pm(k) \calM(N, k) = \calM(N, k + (\pm 2, \mp 1)) C(N, k + (\pm 2, \mp 1))^{-1} U^{-1} H(N, \widehat{E}_\pm(k)) UC(N, k),
\end{equation*}
for all $N \geq 0$. Writing $c = c_N(k)$ and $c_\pm = c_N(k + (\pm 2, \mp 1))$, we have that
\begin{align*}
H(N, \widehat{E}^U_+(k)) &= \begin{pmatrix} N + k_2 - \frac{1}{2} & c(N + k_2 - \frac{1}{2}) + 1 - c_+(N + k_2 + \frac{1}{2}) \\ 0 & N + k_2 + \frac{1}{2} \end{pmatrix}, \\
H(N, \widehat{E}^U_-(k)) &= \begin{pmatrix} N + k_1 + k_2 - \frac{1}{2} & c(N + k_1 + k_2 - \frac{1}{2}) - 1 - c_-(N + k_1 + k_2 + \frac{1}{2}) \\ 0 & N + k_1 + k_2 + \frac{1}{2} \end{pmatrix}.
\end{align*}
Using that
\begin{equation*}
c_\pm = \frac{k_1 \pm 2}{1 + 2N + k_1 + 2k_2},
\end{equation*}
we find that
\begin{align*}
c(N + k_2 - \tfrac{1}{2}) + 1 - c_+(N + k_2 + \tfrac{1}{2}) &= 0, \\
c(N + k_1 + k_2 - \tfrac{1}{2}) - 1 - c_-(N + k_1 + k_2 + \tfrac{1}{2}) &= 0,
\end{align*}
which completes the proof.
\end{proof}


\subsection{General shift operators for the $\bbM_2$-valued Jacobi polynomials}\label{S:gen_mv_shiftops}
Now that we have matrix-valued analogs of the fundamental symmetric shift operators, we can construct matrix-valued shift operators with arbitrary shifts $\ell \in 2\bbZ \times \bbZ$. Multiplying these operators from the right with Cherednik operators also produces matrix-valued shift operators. In fact, we show that any matrix-valued shift operator arises this way, which motivates the terminology of fundamental matrix-valued shift operators.

\begin{lemma}\label{lem:general_G_hat_explicit}
Let $\ell \in 2\bbZ \times \bbZ$ and define the operator $\widehat{G}(\ell, k)$ as in Proposition \ref{prop:general_g} in term of $\widehat{G}_\pm(k)$ and $\widehat{E}_\pm(k)$. Then $\widehat{G}(\ell, k)$ is a $\frakD_k$-shift operator, which is given explicitly by
\begin{equation*}
\widehat{G}(\ell, k) =
\begin{cases}
\begin{pmatrix}
G(\ell, k_+) & 0 \\
-\frac{1}{2}\ell_1 G(\ell - (2, -1), k_+) & G(\ell, k_-)
\end{pmatrix} & \text{if } \ell_1 \geq 0, \\
\begin{pmatrix}
G(\ell, k_+) & \frac{1}{2}\ell_1 G(\ell + (2, -1), k_-) \\
0 & G(\ell, k_-)
\end{pmatrix} & \text{if } \ell_1 < 0.
\end{cases}
\end{equation*}
\end{lemma}

\begin{proof}
This follows from Proposition \ref{prop:frakD_shift_ops}.
\end{proof}

The following lemma allows us to state the analog of Theorem \ref{thm:sym_shift_str_thm}(ii) for $\calD_k$- (or $\frakD_k$-) shift operators in terms of the operator $\calD_k$ (or $\frakD_k$).

\begin{lemma}\label{lem:diag_MLk_as_poly}
The differential operators $ML_{k_\pm}$ can be recovered from $\frakD_k$, i.e.,
\begin{equation*}
\diag(ML_{k_+}, ML_{k_-}) = \frakD_k^2 - \frakD_k.
\end{equation*}
\end{lemma}

\begin{proof}
Follows from a simple computation using \eqref{BC1:frakD}.
\end{proof}

We prove the matrix-valued analogs of Theorem \ref{thm:sym_shift_str_thm}(ii) by reducing it to the scalar-case.

\begin{theorem}\label{thm:mv_structure_thm}
Any $\frakD_k$-shift operator with shift $\ell$ is of the form
\begin{equation*}
\widehat{G}(\ell, k) p(\calD_k), \quad p \in \bbC[\xi].
\end{equation*}
\end{theorem}

\begin{proof}
Let $D$ be a $\frakD_k$-shift operator with shift $\ell$. We only prove the case of $\ell \in 2\bbZ_{\geq 0} \times \bbZ$, as the proof of the other case is analogous. Comparing the entries of the transmutation property yields
\begin{align*}
(1, 1): \quad D_{12}E_-(k_+) &= E_+(k_- + \ell)D_{21} + \tfrac{1}{2}\ell_1 D_{11}, \\
(1, 2): \quad D_{11}E_+(k_-) &= E_+(k_- + \ell)D_{22} + (k_1 + \ell_1) D_{12}, \\
(2, 1): \quad D_{22}E_-(k_+) &= E_-(k_+ + \ell)D_{11} - (k_1 + \ell_1) D_{21}, \\
(2, 2): \quad D_{21}E_+(k_-) &= E_-(k_+ + \ell)D_{12} - \tfrac{1}{2}\ell_1 D_{22}.
\end{align*}
These four equations imply that the entries of $D$ are symmetric shift operators, i.e.,
\begin{equation}\label{eq:four_transmutations}
\begin{aligned}
D_{11}ML_{k_+} &= ML_{k_+ + \ell} D_{11}, \\
D_{12}ML_{k_-} &= ML_{k_- + (2, -1) + \ell} D_{12}, \\
D_{21}ML_{k_+} &= ML_{k_+ - (2, -1) + \ell} D_{21}, \\
D_{22}ML_{k_-} &= ML_{k_- + \ell} D_{22}.
\end{aligned}
\end{equation}
By Theorem \ref{thm:sym_shift_str_thm} there are polynomials $p_{11}$, $p_{12}$, $p_{21}$, and $p_{22}$ such that
\begin{equation}\label{eq:shift_ops}
\begin{aligned}
D_{11} &= G(\ell, k_+) p_{11}(ML_{k_+}), \\
D_{12} &= G(\ell + (2, -1), k_-) p_{12}(ML_{k_-}), \\
D_{21} &= G(\ell - (2, -1), k_+) p_{21}(ML_{k_+}), \\
D_{22} &= G(\ell, k_-) p_{22}(ML_{k_-}).
\end{aligned}
\end{equation}

Let $N = \ord(D)$ and $N' = \ord(G(\ell, k_\pm)) = \ord(G(\ell \pm (2, -1), k_\mp)) \mp 1$. Suppose that $N - N' = 2M$ is even. Since $\ord(ML_k) = 2$, we find that $\max(\ord(D_{12}), \ord(D_{21})) \leq N - 1$ and $\max(\ord(D_{11}), \ord(D_{22})) = N$. From $(1, 2)$ we see that $\ord(D_{11}) = \ord(D_{22}) = N$. Thus $\deg(p_{11}) = \deg(p_{22}) = M$. Moreover, using $(1, 2)$ we see that the leading coefficients of $D_{11}$ and $D_{22}$ coincide, which we denote by $c$. If $M = 0$ then $D_{12} = 0$ because $\ord(G(\ell + (2, -1), k_-)) = N' + 1 > N$. Similarly, it follows that $\deg(p_{12}) = 0$, and by $(1, 1)$ we have that the leading coefficients of $D_{21}$ equals $\frac{1}{2}\ell_1 c$. Hence, $p_{11} = p_{22} = c$, $p_{12} = -\frac{1}{2}\ell_1 c$, and $p_{21} = 0$. Using Lemma \ref{lem:general_G_hat_explicit} we see that
\begin{equation*}
D = c\widehat{G}(\ell, k).
\end{equation*}
Otherwise $M \geq 1$ and by Lemma \ref{lem:general_G_hat_explicit} we have that
\begin{equation*}
D' = D - c\widehat{G}(\ell, k) \diag(ML_{k_+}, ML_{k_-})^M
\end{equation*}
is a $\frakD_k$-shift operator with shift $\ell$, and that $\ord(D') \leq N - 1$, so we can proceed by induction. Note that Lemma \ref{lem:diag_MLk_as_poly} guarantees us that $\diag(ML_{k_+}, ML_{k_-})^M$ can be expressed as a polynomial in $\calD_k$.

Suppose that $N - N' = 2M + 1$ is odd. In this case we find that $\max(\ord(D_{11}), \ord(D_{22}))$ $\leq N - 1$ and $\ord(D_{12}) = \ord(D_{21}) = N$. Thus $\deg(p_{12}) = M$ and $ \deg(p_{21}) = M + 1$. By $(1, 1)$ the leading coefficients of $D_{12}$ and $D_{21}$ both equal some constant $c$. Using \eqref{BC1:frakD} and Lemma \ref{lem:general_G_hat_explicit}, we compute
\begin{align*}
\widehat{G}(\ell, k) \calD_k &=
\begin{pmatrix}
G(\ell, k_+) & -G(\ell, k_+) E_+(k_-) \\
-G(\ell, k_-)E_-(k_+) & G(\ell, k_-)
\end{pmatrix}
+ \text{lower order terms} \\
&=
\begin{pmatrix}
G(\ell, k_+) & -G(\ell + (2, -1), k_-) \\
-G(\ell - (2, -1), k_+) ML_{k_+} & G(\ell, k_-)
\end{pmatrix}
+ \text{lower order terms},
\end{align*}
because $G(\ell, k_+) E_+(k_-) = G(\ell + (2, -1), k_-)$ and 
\begin{align*}
G(\ell, k_-)E_-(k_+) &= G(\ell - (2, -1), k_+) E_+(k_-) E_-(k_+) \\
&= G(\ell - (2, -1), k_+) ML_{k_+} + \text{lower order terms}.
\end{align*}
It follows that 
\begin{equation*}
D' = D + c\widehat{G}(\ell, k) \calD_k \diag(ML_{k_+}, ML_{k_-})^M
\end{equation*}
is a $\frakD_k$-shift operator with shift $\ell$ and $\ord(D') \leq N - 1$, so we can proceed by induction once more.
\end{proof}

\begin{corollary}[of Theorem \ref{thm:mv_structure_thm}]\label{cor:mv_structure_thm_u}
Any $\calD_k$-shift operator with shift $\ell$ is of the form
\begin{equation*}
\widehat{G}^U(\ell, k) p(\frakD_k), \quad p \in \bbC[\xi].
\end{equation*}
\end{corollary}

\begin{corollary}\label{cor:diag_shift_factors}
Any $\calD_k$-shift operator with shift $\ell$ has diagonal shift factors.
\end{corollary}

\begin{proof}
This follows from \eqref{eq:BC1DIffOpM}, Proposition \ref{prop:matrix_G_shift_E}, and Corollary \ref{cor:mv_structure_thm_u}.
\end{proof}


\subsection{The algebra $\calD(\calW_k)$}\label{subsection:DW}
The previous discussion can be extended to determine the algebra of differential operators $\calD(\calW_k)$ associated with a family of matrix-valued orthogonal polynomials for the weight $\calW_k$. More precisely, let $\bbA_2 = \bbC[x, \partial_x] \otimes \bbM_2$ be the Weyl algebra acting on $\bbM_2[x]$, and let $\{ \calQ_n \}_{n \geq 0}$ be any family of matrix-valued orthogonal polynomials for the matrix-weight $\calW_k$. Then
\begin{equation*}
\calD(\calW_k) = \{ D \in \bbA_2 \mid D\calQ_n = \calQ_n \Lambda_n(D), \text{ with } \Lambda_n(D) \in \bbM_2, \text{ for all } n \geq 0 \}.
\end{equation*} 
The algebra $\calD(\calW_k)$ is independent of the particular family $\{ \calQ_n \}_{n \geq 0}$.

The elements
\begin{equation*}
I_2 = \begin{pmatrix} 1 & 0 \\ 0 & 1 \end{pmatrix}, \quad S = \begin{pmatrix} 0 & 1 \\ 1 & 0 \end{pmatrix}, \quad \calD_k
\end{equation*}
are contained in $\calD(\calW_k)$ because $S' = \diag(1,-1) \in \calD(U\calW_k U^*) = \calD(\diag(w_{k_+}, w_{k_-}))$. Note that $\calD(\calW_k)$ is non-commutative since $S\calD_k = -\calD_kS + k_1 I_2 + S \neq \calD_k S$, and that $S \in \calD(\calW_k)$ implies two identities for the non-symmetric Jacobi polynomials.

\begin{corollary}
For $N \geq 0$, we have
\begin{align*}
E(-N,k)(z^{-1}) &= z^{-1}(E(N+1,k)(z)-c_{N}(k)E(-N,k)(z)), \\
E(N+1,k)(z^{-1}) &= z^{-1}(c_{N}(k)E(N+1,k)(z)+(1-c_{N}(k)^2)E(-N,k)(z)).
\end{align*}
\end{corollary}

\begin{proof}
Note that $S' \calN(N, k) = \calN(N, k) S'$, and combining this with Proposition \ref{prop:E_P_rel} gives
\begin{equation*}
S \calM(N, k) = \calM(N, k) C(N, k)^{-1} S C(N, k),
\end{equation*}
which implies
\begin{equation*}
\calM(N, k) = \begin{pmatrix} 0 & 1 \\ 1 & 0 \end{pmatrix} \calM(N, k) \begin{pmatrix}-c_N(k) & 1 - c_N(k)^2 \\ 1 & c_N(k) \end{pmatrix}.
\end{equation*}
Looking at the entries, we find
\begin{align*}
\underline{E(-N, k)}_1 &= \underline{E(N + 1, k)}_2 - c_N(k) \underline{E(-N, k)}_2, \\
\underline{E(-N, k)}_2 &= \underline{E(N + 1, k)}_1 - c_N(k) \underline{E(-N, k)}_1, \\
\underline{E(N + 1, k)}_1 &= c_N(k)\underline{E(N + 1, k)}_2 + (1 - c_N(k)^2) \underline{E(-N, k)}_2, \\
\underline{E(N + 1, k)}_2 &= c_N(k)\underline{E(N + 1, k)}_1 + (1 - c_N(k)^2) \underline{E(-N, k)}_1,
\end{align*}
which imply the desired identities for $E(n, k)(z^{-1}) = \underline{E(n, k)}_1 + \underline{E(n, k)}_2z^{-1}$, with $n \in \bbZ$, since $E(n, k)(z) = \underline{E(n, k)}_1 + \underline{E(n, k)}_2z$, see Lemma \ref{lem:steinberg_basis}.
\end{proof}

We modify the proof of Theorem \ref{thm:mv_structure_thm} to determine $\calD(\diag(w_{k_+}, w_{k_-}))$.

\begin{proposition}\label{prop:DWdiag}
The algebra $\calD(\diag(w_{k_+}, w_{k_-}))$ is generated by $I$, $S'$, and $\frakD_k$.
\end{proposition}

\begin{proof}
Let $D \in \calD(\diag(w_{k_+}, w_{k_-}))$ be a differential operator of order $N$. The relation $D\calN(N',k) = \calN(N',k) H(N', D)$ gives
\begin{align*}
(1, 1): \quad D_{11}p(N', k_+) &= H(N', D)_{11} p(N', k_+), \\
(1, 2): \quad D_{12}p(N', k_-) &= H(N', D)_{12} p(N', k_+), \\
(2, 1): \quad D_{21}p(N', k_+) &= H(N', D)_{21} p(N', k_-), \\
(2, 2): \quad D_{22}p(N', k_-) &= H(N', D)_{22} p(N', k_-)
\end{align*}
for all $N' \geq 0$. It follows that \eqref{eq:four_transmutations} holds with $\ell = 0$ for $D$. Thus \eqref{eq:shift_ops} also holds with $\ell = 0$.

Suppose that $N = 2M$ is even. As before, we have $\max(\ord(D_{12}), \ord(D_{21})) \leq N - 1$ and $\max(\ord(D_{11}), \ord(D_{22})) = N$. However, it may happen that $\ord(D_{11}) \neq \ord(D_{22})$. If $M = 0$ then $D_{12} = D_{21} = 0$, so that $D$ is a diagonal constant matrix. Clearly, $D$ is a linear combination of $I$ and $S'$. Otherwise $M \geq 1$ and there exist constants $c$ and $c'$ such that $c + c'$ and $c - c'$ equal the leading coefficients of $D_{11}$ and $D_{22}$, respectively. Thus
\begin{equation*}
D' = D - c \thinspace \diag(ML_{k_+}, ML_{k_-})^M - c'S' \thinspace \diag(ML_{k_+}, ML_{k_-})^M
\end{equation*}
is an element of $\calD(\calW_k)$ and $\ord(D') \leq N - 1$, so we can proceed by induction. By Lemma \ref{lem:diag_MLk_as_poly} we have that $D' - D$ lies in the algebra generated by $I$, $S'$, and $\calD_k$.

Now suppose that $N = 2M + 1$ is odd. In this case we have $\max(\ord(D_{11}), \ord(D_{22})) \leq N - 1$, $\max(\ord(D_{12}), \ord(D_{21})) = N$, and possibly $\ord(D_{12}) \neq \ord(D_{21})$. Again, there exist constants $c$ and $c'$ such that $c + c'$ and $c - c'$ equal the leading coefficients of $D_{12}$ and $D_{21}$, respectively. Hence
\begin{equation*}
D' = D + \calD_k (c \thinspace \diag(ML_{k_+}, ML_{k_-})^M + c'S' \thinspace \diag(ML_{k_+}, ML_{k_-})^M)
\end{equation*}
is an element of $\calD(\calW_k)$ and $\ord(D') \leq N - 1$, so we can proceed by induction one last time. Note that $D - D'$ is an element of the algebra generated by $I$, $S'$, and $\calD_k$.
\end{proof}

Let us state some of the corollaries of (the proof of) Proposition \ref{prop:DWdiag}.

\begin{corollary}
The commutant of $\diag(ML_{k_+}, ML_{k_-})$ in the Weyl algebra $\bbA_2$ equals the algebra $\calD(\diag(w_{k_+}, w_{k_-}))$.
\end{corollary}

The algebra $\bbA_2$ is filtered by degree, i.e., 
\begin{equation*}
\bbA_2 = \bigcup_{N\ge0}\bbA_2^{(N)}, \quad \bbA_2^{(N)}=\{D \in \bbA_2 \mid \ord(D) \leq N \},
\end{equation*}
and this induces a filtration on the subalgebra $\calD(\calW_k)$.

\begin{corollary}\label{cor:DW}
The algebra $\calD(\calW_k)$ is generated by $I$, $S$, and $\calD_k$. Moreover, each quotient $\calD(\calW_k)^{(N)}/\calD(\calW_k)^{(N - 1)}$ is a two-dimensional vector space.
\end{corollary}

\begin{remark}
In the Chebyshev case, i.e., $k = 0$, the matrix-valued orthogonal polynomials $\calM(N, 0)$ have been introduced in \cite{Ber68}, and studied further in \cite{CaGr05, CaGr06}. The explicit expression for $\calM(N, 0)$ from \cite[Prop.2.5]{vHvP1} reduces to \cite[(3.14)]{Ber68}, see also \cite[(8.4)]{CaGr06}. The operator $\calD_0$ is closely related to the operator appearing in \cite[(3.2)]{CaGr05}, and $S$ is also given in \cite[(3.3)]{CaGr05}. Finally, the algebra $\calD(\calW_0)$ has been determined, albeit without proof, in \cite[Third example]{CaGr06}, and Corollary \ref{cor:DW} can be viewed as an extension of this result.
\end{remark}


\subsection{Shift operators for the $\bbC^2$-valued Jacobi polynomials}
Finally, we interpret each $\calD_k$-shift operator as an operator on $\bbC[x] \otimes \bbC^2$, and Corollary \ref{cor:diag_shift_factors} allows us to describe the action on the vector-valued Jacobi polynomials.

\begin{corollary}\label{cor:calD_shifts_act_on_vvjp}
Viewing any $\calD_k$-shift operator $\widehat{S}^U(k)$ with shift $\ell$ as a differential operator on $\bbC[x] \otimes \bbC^2$, we have for $N \geq 0$ that
\begin{align*}
\widehat{S}^U(k) \calP(-N,k) &= H(N, \widehat{S}^U(k))_{11} \thinspace \calP(-(N - \rho(\ell)), k), \\
\widehat{S}^U(k) \calP(N+1,k) &= H(N, \widehat{S}^U(k))_{22} \thinspace \calP(N + 1 - \rho(\ell), k).
\end{align*}
\end{corollary}

Recall from \eqref{eq:mv_inner_prod_prime} that
\begin{equation}\label{eq:mv_inprod_rel}
(U^{-1}\calP, U^{-1}\calQ)_k = (\calP, \calQ)'_k = (\calP_1, \calQ_1)_{k_+} + (\calP_2, \calQ_2)_{k_-},
\end{equation}
and note that on the right-hand side we have two scalar-valued inner products, so that the formal adjoints of the $\frakD_k$-shift operators can be computed from those of their entries.

\begin{lemma}\label{lem:matrix_valued_adjoints}
For $\ell \in 2\bbZ \times \bbZ$, we have
\begin{equation*}
(\widehat{G}(\ell, k)\calP, \calQ)'_{k + \ell} = (\calP, \widehat{G}(-\ell, k + \ell)\calQ)'_k, \quad \calP, \calQ \in \bbC[x] \otimes \bbC^2.
\end{equation*}
\end{lemma}

\begin{proofempty}
We present the proof for $\ell \geq 0$. Using \eqref{eq:mv_inprod_rel} and Lemma \ref{lem:sym_g_adjoint}, we compute
\pushQED{\qed}
\begin{align*}
&(\widehat{G}(\ell, k)\calP, \calQ)'_{k + (2, -1)} \\
&= (G(\ell, k_+)\calP_1, \calQ_1)_{k_+ + \ell} + (G(\ell, k_-)\calP_2 - \tfrac{1}{2} \ell_1 G(\ell - (2, -1), k_+) \calP_1, \calQ_2)_{k_- + \ell} \\
&= (\calP_1, G(-\ell, k_+ + \ell)\calQ_1 - \tfrac{1}{2} \ell_1 G(-\ell + (2, -1), k_- + \ell) \calQ_2)_{k_+} + (\calP_2, G(-\ell, k_- + \ell)\calQ_2)_{k_-} \\
&= (\calP, \widehat{G}(-\ell, k + \ell)\calQ)'_k. \qedhere
\end{align*}
\popQED
\end{proofempty}

\begin{corollary}\label{cor:formal_adj_g_hat}
If $\widehat{S}^U(k) = \widehat{G}^U(\ell, k)p(D_k)$ is any $\calD_k$-shift operator with shift $\ell$, then
\begin{equation*}
(\widehat{G}^U(\ell, k)p(D_k)\calP, \calQ)_{k + \ell} = (\calP, \widehat{G}^U(-\ell, k + \ell)p(D_{k + \ell})\calQ)_k, \quad \calP, \calQ \in \bbC[x] \otimes \bbC^2.
\end{equation*}
\end{corollary}

\begin{proof}
This follows from \eqref{eq:mv_inprod_rel}, Lemma \ref{lem:matrix_valued_adjoints}, and that of $D_k$ is formally self-adjoint.
\end{proof}


\section{Differential-reflection shift operators}\label{S:ShiftNonSym}


\subsection{Fundamental shift operators for non-symmetric Jacobi polynomials}
Any $\calD_k$-shift operator acts on the vector-valued Jacobi polynomials by Corollary \ref{cor:calD_shifts_act_on_vvjp}. We consider the pullback of a $\calD_k$-shift operator through the map $\Upsilon : \bbC[z^{\pm 1}] \to \bbC[x] \otimes \bbC^2$, which acts on the non-symmetric Jacobi polynomials by construction. For this we introduce the notion of a non-symmetric shift operator.

\begin{definition}\label{def:nonsym_shift_op}
An element $\calS(k) \in \bbC_\Delta[z^\pm] \otimes \bbC[z\partial_z] \otimes \bbC[\bbZ_2]$ is called a non-symmetric shift operator with shift $\ell \in 2\bbZ \times \bbZ$ if
\begin{enumerate}[(i)]
    \item $\calS(k)$ satisfies the transmutation property, i.e., $\calS(k) D_k = D_{k + \ell} \calS(k)$,
    \item $\calS(k)$ acts as a differential-reflection operator on $\bbC[z^{\pm 1}]$.
\end{enumerate}
\end{definition}

\begin{lemma}\label{lem:pullback_shifts}
There is a $1$-to-$1$ correspondence between non-symmetric shift operators with shift $\ell$ and $\calD_k$-shift operators with shift $\ell$, given by the pullback and pushforward through $\Upsilon$.
\end{lemma}

\begin{proof}
This follows from the fact that $f_1$ and $f_2$ can be expressed in terms of $f$ using the reflection $s$, see Lemma \ref{lem:steinberg_basis}.
\end{proof}

\begin{proposition}\label{prop:nonsym_shifts_expl}
Let $\calG_\pm(k)$ and $\calE_\pm(k)$ be the $\Upsilon$-pullbacks of $\widehat{G}^U_\pm(k)$ and $\widehat{E}^U_\pm(k)$. Explicitly, the non-symmetric forward and backward shift operators are given by
\begin{equation*}
\calG_+(k) = G_+(k) - \frac{z(1 - s)}{(z - z^{-1})^2}, \quad \calG_-(k) = G_-(k) + z^{-1} - z s,
\end{equation*}
and the non-symmetric contiguity shift operators by
\begin{equation*}
\calE_+(k) = E_+(k) - \frac{1 + z^{-1}}{1 - z^{-1}} \frac{1 - s}{z - z^{-1}}, \quad \calE_-(k) = E_-(k) + \frac{1 - z^{-1}}{1 + z^{-1}} \frac{1 - s}{z - z^{-1}}.
\end{equation*}
These four operators are called the fundamental non-symmetric shift operators. The shifts of $\calG_\pm(k)$ and $\calE_\pm(k)$ are given by $(0, \pm 1)$ and $(\pm 2, \mp 1)$, respectively.
\end{proposition}

\begin{proof}
It follows from Lemma \ref{lem:pullback_shifts} that $\calG_\pm(k)$ and $\calE_\pm(k)$ are non-symmetric shift operators with shifts $(0, \pm 1)$ and $(\pm 2, \mp 1)$, respectively. The explicit expressions are derived in Appendix \ref{app:a}.
\end{proof}

\begin{remark}
Let $\calS(k)$ be a non-symmetric shift operator with shift $\ell$. Writing $\mathrm{Sym}(\calS(k)) = (\sum_{w \in \bbZ_2} w\calS(k))$ we have $\mathrm{Sym}(\calS(k)) D_k^2 = D_{k + \ell}^2 \mathrm{Sym}(\calS(k))$. Since $\mathrm{Sym}(\calS(k))$ restricts to an operator on $\bbC[x]$, we can apply $\beta$ to obtain
\begin{equation*}
\beta(\mathrm{Sym}(\calS(k)))ML_k = ML_{k + \ell} \beta(\mathrm{Sym}(\calS(k))).
\end{equation*}
It follows that $\beta(\mathrm{Sym}(\calS(k)))$ is a symmetric shift operator. This can also be verified directly from Proposition \ref{prop:nonsym_shifts_expl}, and the fundamental symmetric shift can be recovered from this.
\end{remark}

To present uniform expressions for the action of the fundamental non-symmetric shift operators on the non-symmetric Jacobi polynomials, we introduce a new labeling for the non-symmetric Jacobi polynomials. We employ the notation of \cite{Opd95} but specialized to the root system of type $BC_1$. Recall that $P = \bbZ \eps$, $P^+ = \bbZ_{\geq 0}\eps$, and $s\eps = -\eps$. The elements of $\bbZ_2$ are ordered by their word length, i.e., $\ell(1) = 0$ and $\ell(s) = 1$. For $\lambda \in P$ we have that its $\bbZ_2$-orbit equals $\{ \lambda, -\lambda \}$ which intersects $P^+$ exactly once. We define the element $\lambda^* \in P^+$ by $\bbZ_2 \cap P^+ = \{ \lambda^* \}$ and $w^\lambda \in \bbZ_2$ to be the shortest element such that $\lambda = w^\lambda \lambda^*$. Finally, we denote the longest of the stabilizer of $\lambda$ in $\bbZ_2$ by $w_\lambda$. Explicitly, for $\lambda_n = n\eps$ we have $\lambda_n^* = |n|\eps$, $w^{\lambda_n} = 1$ iff $n \in \bbZ_{\geq 0}$, and $w_{\lambda_n} = s$ iff $n = 0$. Finally, we define
\begin{equation*}
E(\lambda_n, k) = E(n, k), \quad \eta(\lambda_n, k) = \eta(n, k).
\end{equation*}

\begin{lemma}\label{lem:nonsym_shift_factors}
If $\calS(k)$ is a fundamental non-symmetric shift operator with shift $\ell$, then
\begin{equation*}
\calS(k)E(\lambda, k) =
\begin{cases}
\eta(\lambda, \calS(k)) E(w^\lambda w_{\lambda^*} (\lambda^* - \rho(\ell)\eps), k + \ell) & \text{if } \lambda^* - \rho(\ell)\eps \in P^+, \\
0 & \text{otherwise},
\end{cases}
\end{equation*}
with shift factors
\begin{equation*}
\begin{aligned}
\eta(-n, \calG_+(k)) &= n, \qquad\qquad\qquad\;\;\; \thinspace \eta(n + 1, \calG_+(k)) = n, \\
\eta(-n, \calG_-(k)) &= n + k_1 + 2k_2, \qquad \thinspace \eta(n + 1, \calG_-(k)) = n + k_1 + 2k_2, \\
\eta(-n, \calE_+(k)) &= n + k_2 - \tfrac{1}{2}, \qquad\quad \thinspace \eta(n + 1, \calE_+(k)) = n + k_2 + \tfrac{1}{2}, \\
\eta(-n, \calE_-(k)) &= n + k_1 + k_2 - \tfrac{1}{2}, \quad \eta(n + 1, \calE_-(k)) = n + k_1 + k_2 + \tfrac{1}{2}.
\end{aligned}
\end{equation*}
\end{lemma}

\begin{proof}
This follows from Proposition \ref{prop:matrix_G_shift_E}.
\end{proof}


\subsection{General shift operators for non-symmetric Jacobi polynomials}
In view of Lemma \ref{lem:pullback_shifts}, we can transfer the results about $\calD_k$-shift operators from Section \ref{S:gen_mv_shiftops} to statements concerning non-symmetric shift operators.

\begin{proposition}\label{prop:pre_nonsym_shift_str_thm}
Let $\ell \in 2\bbZ \times \bbZ$ and define the operator $\calG(\ell, k)$ as in Proposition \ref{prop:general_g} in term of $\calG_\pm(k)$ and $\calE_\pm(k)$, or equivalently as the $\Upsilon$-pullback of $\widehat{G}^U_+(\ell, k)$. Then $\calG(\ell, k)$ is a non-symmetric shift operator, and any element of $\calS(k) \in \bbS_\ns(\ell, k)$ is of the form
\begin{equation*}
\calS(k) = \calG(\ell, k) p(D_k), \quad p \in \bbC[\xi].
\end{equation*}
\end{proposition}

\begin{proof}
This follows from Lemma \ref{lem:general_G_hat_explicit} and Corollary \ref{cor:mv_structure_thm_u}.
\end{proof}

\begin{corollary}\label{cor:nonsym_shiftop_action_and_adjoint}
Let $\calS(k) = \calG(\ell, k)p(D_k)$ be any non-symmetric shift operator with shift $\ell$.
\begin{enumerate}[(i)]
    \item The action of $\calS(k)$ on the non-symmetric Jacobi polynomials is given by
    \begin{equation*}
    \calS(k)E(\lambda, k) =
    \begin{cases}
    \eta(\lambda, \calS(k)) E(w^\lambda w_{\lambda^*} (\lambda^* - \rho(\ell)\eps), k + \ell) & \text{if } \lambda^* - \rho(\ell)\eps \in P^+, \\
    0 & \text{otherwise},
    \end{cases}
    \end{equation*}
    with shift factors $\eta(\lambda, \calS(k)) = \eta(\lambda, \calG(\ell, k))p(\widetilde{\lambda})$. Here $\widetilde{n\eps} = n + \rho(k)$ for $n > 0$ and $\widetilde{n\eps} = n - \rho(k)$ and $n \leq 0$.
    \item The formal adjoint of $\calS(k)$ is given by $\calG(-\ell, k + \ell)p(D_{k + \ell})$, i.e.,
    \begin{equation*}
    (\calG(\ell, k)p(D_k)f, g)_{k + \ell} = (f, \calG(-\ell, k + \ell)p(D_{k + \ell})g)_k, \quad f, g \in \bbC[z^{\pm 1}].
    \end{equation*}
\end{enumerate}
\end{corollary}

\begin{proof}
Statement (i) follows from \eqref{eq:BC1_eigenvalue} and Corollary \ref{cor:calD_shifts_act_on_vvjp}, and (ii) from Corollary \ref{cor:formal_adj_g_hat}.
\end{proof}

Let $\bbS_\ns(\ell, k)$ denote space of non-symmetric shift operators with shift $\ell \in 2\bbZ \times \bbZ$. The composition of $\calS(k) \in \bbS_\ns(\ell, k)$ and $\calS'(k + \ell) \in \bbS_\ns(\ell', k + \ell)$ gives $\calS'(k + \ell) \calS(k) \in \bbS_\ns(\ell + \ell', k)$. This implies that $\bbS_\ns(\ell, k)$ is a right $\bbS_\ns(0, k)$-module, and that $\bbS_\ns(0, k)$ is an algebra, cf. Section \ref{S:shiftops_for_sym_jac}. Proposition \ref{prop:pre_nonsym_shift_str_thm} has interesting consequences concerning the structure of $\bbS_\ns(\ell, k)$. We summarize these corollaries into the so-called structure theorem for the non-symmetric shift operators, cf. Theorem \ref{thm:sym_shift_str_thm}.

\begin{theorem}[Structure theorem for $\bbS_\ns(\ell, k)$]\label{thm:nonsym_shift_str_thm}\text{ }
\begin{enumerate}[(i)]
    \item The algebra $\bbS_\ns(0, k)$ is generated by the Cherednik operator $D_k$.
    \item The $\bbS_\ns(0, k)$-module $\bbS_\ns(\ell, k)$ is of rank one and freely generated by $\calG(\ell, k)$.
    \item The generators satisfy the relations
    \begin{equation*}
    \calG(\ell + \ell', k) = \calG(\ell', k + \ell) \calG(\ell, k) = \calG(\ell, k + \ell') \calG(\ell', k).
    \end{equation*}
    whenever $\sign(\ell_1) = \sign(\ell'_1)$ and $\sign(\frac{1}{2}\ell_1 + \ell_2) = \sign(\frac{1}{2}\ell'_1 + \ell'_2)$. 
\end{enumerate}
\end{theorem}


\subsection{Alternative proof of Proposition \ref{prop:pre_nonsym_shift_str_thm} via differential-reflection operators}\label{S:SpaceOfShiftOps}
We present another proof of Proposition \ref{prop:pre_nonsym_shift_str_thm} that is more fitting in the context of differential-reflection operators, and it is the natural adaption of the proof of \cite[Prop.3.3.4]{HeSc94}. This proof is independent of the structure theorem for the symmetric shift operators. We deduce from Definition \ref{def:nonsym_shift_op} that each element of $\bbS_\ns(\ell, k)$ admits an order, and that $\calG(\ell, k)$ is the smallest element, up to a scalar multiple, with respect to this order. This induces a filtration on $\bbS_\ns(\ell, k)$ that can be used to prove Proposition \ref{prop:pre_nonsym_shift_str_thm} inductively.

\begin{lemma}\label{lem:tilde_nonsym_shifts}
For $\calS(k) \in \bbS_\ns(\ell, k)$, we have $\widetilde{\calS}(k) = \delta_{-k} \circ \calS(k)^* \circ \delta_{k + \ell} \in \bbS_\ns(-\ell, k + \ell)$, and
\begin{equation*}
(\calS(k)f, g)_{k + \ell} = (f, \widetilde{\calS}(k)g)_k, \quad f,g \in \bbC[z^{\pm 1}].
\end{equation*}
Here $*$ denotes the formal adjoint of a differential-reflection operator on $S^1$ with respect \eqref{eq:pairing_on_s1}.
\end{lemma}

\begin{proof}
The proof of \cite[Prop.3.4]{Opd88a} only relies on the fact that $ML_k$ is formally self-adjoint, which is also the case for $D_k$, so it also applies to this case.
\end{proof}

\begin{proposition}\label{prop:order_of_nonsym_shiftop}
For $\calS(k) \in \bbS_\ns(\ell, k) \backslash \{ 0 \}$, there exist unique $N \geq 0$ and $a(k) \in \bbC \backslash \{ 0 \}$ such that
\begin{equation*}
\calS(k) = a(k)(z^{\tfrac{1}{2}} - z^{-\tfrac{1}{2}})^{-\ell_1} (z - z^{-1})^{-\ell_2} (z\partial_z)^N + \sum_{i \leq N - 1, \; j = 0, 1} c_{ij}(k) (z\partial_z)^i s^j, \quad c_{ij}(k) \in \bbC_\Delta[z^\pm].
\end{equation*}
The number $N = \ord(\calS(k))$ is called the order of $\calS(k)$ and satisfies the inequality
\begin{equation*}
\ord(\calS(k)) \geq |\tfrac{1}{2}\ell_1| + |\tfrac{1}{2}\ell_1 + \ell_2| = \ord(\calG(\ell, k)),
\end{equation*}
with equality if and only if $\calS(k) = a(k)\calG(\ell, k)$.
\end{proposition}

\begin{proof}
First, writing $\theta = z \partial_z$ the transmutation property becomes
\begin{equation*}
[D_k, \calS(k)] = -\left(\tfrac{1}{2} \ell_1 \frac{1 + z^{-1}}{1 - z^{-1}} + \ell_2 \frac{1 + z^{-2}}{1 - z^{-2}} \right)(1 - s) \calS(k).
\end{equation*}
Taking $N, M \geq 0$ and $a_i(k), b_j(k) \in \bbC_\Delta[z^\pm]$ with $a_N(k), b_M(k) \neq 0$ such that
\begin{equation*}
\calS(k) = \sum_{i=0}^N a_i(k) \theta^i + \sum_{j=0}^M b_j(k) \theta^j s,
\end{equation*}
we compute
\begin{equation*}
[D_k, \calS(k)] = \theta(a_N(k)) \theta^N + 2b_M(k) \theta^{M+1} + \sum_{i = 0}^{\max(N - 1, M)} c_i(k) \theta^i + \sum_{j = 0}^{\max(N, M)} d_j(k) \theta^j s,
\end{equation*}
for some $c_i(k), d_j(k) \in \bbC_\Delta[z^\pm]$. Suppose that $M \geq N$ then it follows that $b_M(k) = 0$, which is a contradiction, so we conclude that $M \leq N - 1$. Comparing the coefficients of $\theta^N$ in both expressions of $[D_k, \calS(k)]$ yields
\begin{equation*}
2\theta(a_N(k)) = -a_N(k) \left( \ell_1 \frac{1 + z^{-1}}{1 - z^{-1}} + 2\ell_2 \frac{1 + z^{-2}}{1 - z^{-2}} \right),
\end{equation*}
which is exactly the first-order differential equation that occurs in the proof of \cite[Prop.3.3.4]{HeSc94}, and the solution is
\begin{equation*}
a_N(k) = a(k)(z^{\frac{1}{2}} - z^{-\frac{1}{2}})^{-\ell_1}(z - z^{-1})^{-\ell_2}, \quad a(k) \in \bbC \backslash \{0\}.
\end{equation*}

Second, it will be convenient to consider the operator
\begin{equation*}
\widehat{S}(k) = \Phi^{-1} \circ \begin{pmatrix} \calS(k) & 0 \\ 0 & s\calS(k)s \end{pmatrix} \circ \Phi, \quad \Phi = \begin{pmatrix} 1 & z \\ 1 & z^{-1} \end{pmatrix},
\end{equation*}
which satisfies $\widehat{S}(k)\underline{f} = \underline{\calS(k)f}$ for all $f \in \bbC[z^{\pm 1}]$ because $(f_1, f_2)^T = \Phi^{-1} (f, s(f))^T$. It follows that $\widehat{S}(k)$ is a $\calD_k$-shift operator, and in particular that $\widehat{S}(k)$ is an element of the Weyl algebra $\bbC[x, \partial_x] \otimes \bbM_2$, see Remark \ref{rem:mv_shiftops_weyl_algebra}. At the same time, if $N - \ell_1 - \ell_2$ is even then we compute
\begin{equation*}
\widehat{S}(k) = a_N(k) I_2 \theta^N + \sum_{i = 0}^{N - 1} B_i(k) \theta^i, \quad B_i(k) \in \bbC_\Delta[z^\pm] \otimes \mathbb{M}_2,
\end{equation*}
and find that $a_N(k) \theta^N$ can be expressed in terms of $\partial_x$ and $x$. As in the proof of \cite[Prop.3.3.4]{HeSc94}, we get
\begin{equation*}
a_N(k) \theta^N = 2^{N - \frac{1}{2}\ell_1 - \ell_2} a(k)(x - 1)^{\frac{1}{2}(N - \ell_1 - \ell_2)}(x + 1)^{\frac{1}{2}(N - \ell_2)}\partial_x^N + \sum_{i = 0}^{N - 1} b_i(k)\partial_x^i,
\end{equation*}
which implies that $N - \ell_1 - \ell_2, N - \ell_2 \in 2 \bbZ_{\geq 0}$. Similarly, if $N - \ell_1 - \ell_2$ is odd then
\begin{equation*}
\widehat{S}(k) = \frac{2a_N(k)}{(x - 1)^{\frac{1}{2}}(x + 1)^{\frac{1}{2}}} \begin{pmatrix} -x & -1 \\ 1 & x \end{pmatrix} \theta^N + \sum_{i = 0}^{N - 1} B_i(k) \theta^i, \quad B_i(k) \in \bbC_\Delta[z^\pm] \otimes \mathbb{M}_2,
\end{equation*}
and
\begin{equation*}
\frac{a_N(k)}{(x - 1)^{\frac{1}{2}}(x + 1)^{\frac{1}{2}}} \theta^N = 2^{N - \frac{1}{2}\ell_1 - \ell_2} a(k)(x - 1)^{\frac{1}{2}(N - \ell_1 - \ell_2 - 1)}(x + 1)^{\frac{1}{2}(N - \ell_2 - 1)}\partial_x^N + \sum_{i = 0}^{N - 1} b_i(k)\partial_x^i,
\end{equation*}
so that $N - \ell_1 - \ell_2, N - \ell_2 \in 2\bbZ_{\geq 0} + 1$. In both cases, we conclude that $N - \ell_1 - \ell_2, N - \ell_2 \in \bbZ_{\geq 0}$. 

Repeating the same arguments for $\widetilde{\calS}(k) \in \bbS_\ns(-\ell, k + \ell)$, see Lemma \ref{lem:tilde_nonsym_shifts}, we obtain that $N + \ell_1 + \ell_2, N + \ell_2 \in \bbZ_{\geq 0}$. Together they imply that $N - \max(|\ell_1 + \ell_2|, |\ell_2|) \in \bbZ_{\geq 0}$. Hence
\begin{equation*}
\ord(\calS(k)) \geq \max(|\ell_1 + \ell_2|, |\ell_2|) = |\tfrac{1}{2}\ell_1| + |\tfrac{1}{2}\ell_1 + \ell_2| = \ord(\calG(\ell, k)).
\end{equation*}
The last equality follows directly from the Proposition \ref{prop:pre_nonsym_shift_str_thm}. As for the last statement, suppose that $\ord(\calS(k)) = |\tfrac{1}{2}\ell_1| + |\tfrac{1}{2}\ell_1 + \ell_2|$ then $\ord(\calS(k) - a(k)\calG(\ell, k)) \leq |\tfrac{1}{2}\ell_1| + |\tfrac{1}{2}\ell_1 + \ell_2| - 1$, which cannot occur by the previous inequality when $\calS(k) - a(k)\calG(\ell, k) \neq 0$.
\end{proof}

\begin{remark}
The matrix-valued context was only needed in the proof of Proposition \ref{prop:order_of_nonsym_shiftop} for the argument with the Weyl algebra, see the second part of Remark \ref{rem:mv_shiftops_weyl_algebra}. Notably, we have not used the fact that the matrix-weight diagonalizes.
\end{remark}

We can reformulate Proposition \ref{prop:pre_nonsym_shift_str_thm} slightly in terms of the order on $\bbS_\ns(\ell, k)$ and prove it using Proposition \ref{prop:order_of_nonsym_shiftop}.

\begin{corollary}\label{cor:pre_nonsym_shift_str_thm}
For $\calS(k) \in \bbS_\ns(\ell, k)$, there exists a unique $p \in \bbC[\xi]$ of degree $\ord(\calS(k)) - \ord(\calG(\ell, k))$ such that
\begin{equation*}
\calS(k) = \calG(\ell, k) p(D_k), 
\end{equation*}
\end{corollary}

\begin{proof}
This follows from Proposition \ref{prop:order_of_nonsym_shiftop} by induction on the order of $\calS(k)$, cf. \cite[Prop.3.8]{Opd88a}.
\end{proof}

Finally, Corollary \ref{cor:nonsym_shiftop_action_and_adjoint}(ii) also becomes a consequence of Proposition \ref{prop:order_of_nonsym_shiftop} in conjunction with Lemma \ref{lem:tilde_nonsym_shifts}, which becomes clear after the following corollary.

\begin{corollary}
For $\ell \in 2\bbZ \times \bbZ$, we have $\widetilde{\calG}(\ell, k) = \calG(-\ell, k + \ell)$.
\end{corollary}

\begin{proof}
The leading term of $\widetilde{\calG}(\ell, k) = \delta_{-k} \circ \calG(\ell, k)^* \circ \delta_{k + \ell}$ is given by
\begin{equation*}
(z^{\tfrac{1}{2}} - z^{-\tfrac{1}{2}})^{-\ell_1} (z - z^{-1})^{-\ell_2} (z\partial_z)^{|\tfrac{1}{2}\ell_1| + |\tfrac{1}{2}\ell_1 + \ell_2|},
\end{equation*}
so the statement follows from Proposition \ref{prop:order_of_nonsym_shiftop} and Corollary \ref{cor:pre_nonsym_shift_str_thm}.
\end{proof}


\section{Normalizations of the non-symmetric Jacobi polynomials via non-symmetric shift operators}\label{S:normalization}

In \cite{Opd89}, the norms and the evaluation at the identity of the symmetric Jacobi polynomials for arbitrary root systems have been calculated by Opdam using symmetric shift operators. Both normalizations are expressed in terms of the generalized Harish-Chandra $c$-functions. Similarly, the normalizations of the non-symmetric Jacobi polynomials for arbitrary root systems have also been expressed in terms of the generalized Harish-Chandra $c$-functions in \cite{Opd95}, which is established by reducing it to the known normalizations.

For the root system of type $BC_1$, we recover the formulas of \cite[Thm.5.3]{Opd95} using non-symmetric shift operators.

\begin{lemma}\label{lem:expl_c_funct_bc1}
The generalized Harish-Chandra $c$-functions from \cite[\S5]{Opd95} specialize for the root system of type $BC_1$ to
\begin{align*}
\widetilde{c}_1(\lambda, k) &= \frac{\Gamma(\lambda(e_1))\Gamma(\lambda(e_1) + \frac{1}{2})}{2^{k_1}\Gamma(\lambda(e_1) + \frac{1}{2}k_1 + \frac{1}{2})\Gamma(\lambda(e_1) + \frac{1}{2}k_1 + k_2)}, \\
\widetilde{c}_s(\lambda, k) &= \frac{\Gamma(\lambda(e_1) + \frac{1}{2})\Gamma(\lambda(e_1) + 1)}{2^{k_1}\Gamma(\lambda(e_1) + \frac{1}{2}k_1 + \frac{1}{2})\Gamma(\lambda(e_1) + \frac{1}{2}k_1 + k_2 + 1)},
\end{align*}
and
\begin{align*}
c^*_1(\lambda, k) &= \frac{\Gamma(-\lambda(e_1) - \frac{1}{2}k_1 + \frac{1}{2})\Gamma(-\lambda(e_1) - \frac{1}{2}k_1 - k_2)}{2^{k_1}\Gamma(-\lambda(e_1))\Gamma(-\lambda(e_1) + \frac{1}{2})}, \\
c^*_s(\lambda, k) &= \frac{\Gamma(-\lambda(e_1) - \frac{1}{2}k_1 + \frac{1}{2})\Gamma(-\lambda(e_1) - \frac{1}{2}k_1 - k_2 + 1)}{2^{k_1}\Gamma(-\lambda(e_1) + \frac{1}{2})\Gamma(-\lambda(e_1) + 1)},
\end{align*}
where $\lambda \in P$ and $e_1 \in \bbR$ such that $\eps(e_1) = 1$.
\end{lemma}

\begin{proof}
See the proof of \cite[Thm.3.4.3]{HeSc94} for the first expression. The other expressions can be obtained similarly.
\end{proof}

It will be convenient to introduce the following rescaling of the operators $\calG(\ell, k)$.

\begin{definition}\label{def:nonsym_g_pm}
Let $B = \{(0, 1), (2, -1)\}$ and let $\ell \in \bbZ_{\geq 0}B$. We define
\begin{equation*}
\calG_+(\ell, k) = 4^{\tfrac{1}{2}\ell_1} \calG(\ell, k), \quad \calG_-(-\ell, k) = 4^{\tfrac{1}{2}\ell_1} \calG(-\ell, k),
\end{equation*}
which are called the non-symmetric raising and lowering shift operators, respectively. Similarly, we define the symmetric raising and lowering shift operators $G_+(\ell, k)$ and $G_-(-\ell, k)$.
\end{definition}


\subsection{Norm formula}\label{sec:norm_formula_bc1}
\noindent Using essentially the proof of \cite[Cor.3.5.3]{HeSc94}, we will see that the calculation of the norms of the polynomials $E(n, k)$ reduces to finding expressions for the shift factors of the operators $\calG_\pm(\ell, k)$. We express these shift factors in terms of the generalized Harish-Chandra $c$-functions.

\begin{notation}
For $k \in 2\bbZ \times \bbZ$ and $c \in \bbZ$, the condition $k \geq c$ denotes $k_1 \geq c$ and $k_2 \geq c$.
\end{notation}

\begin{proposition}\label{prop:shift_values_and_c_funct}
Let $\ell \in \bbZ_{\geq 0} B \backslash \{ 0 \}$, $k \geq 1$, and $k - \ell \geq 0$. Then
\begin{align*}
\eta(\lambda, \calG_+(\ell, k)) &=
\frac{c^*_{w^\lambda w_{\lambda^*}}(-(\lambda^* + \rho(k)\eps), k)}{c^*_{w^\lambda w_{\lambda^*}}(-(\lambda^* + \rho(k)\eps), k + \ell)}, \quad \lambda^* - \rho(\ell)\eps \in P^+, \\
\eta(\lambda, \calG_-(-\ell, k)) &= \frac{\widetilde{c}_{w^\lambda w_{\lambda^*}}(\lambda^* + \rho(k)\eps, k - \ell)}{\widetilde{c}_{w^\lambda w_{\lambda^*}}(\lambda^* + \rho(k)\eps, k)}, \qquad \;\; \lambda \in P.
\end{align*}
Moreover, the quotients of generalized Harish-Chandra $c$-functions on the right-hand sides have no poles and no zeros for $\lambda \in P$ satisfying the stated conditions.
\end{proposition}

\begin{proof}
Note that the comment about there not being any poles and zeros follows directly from the expressions of the generalized Harish-Chandra $c$-functions from Lemma \ref{lem:expl_c_funct_bc1}.

We prove the relations for the fundamental (raising and lowering) non-symmetric shift operators. We compute the right-hand side expressions using Lemma \ref{lem:expl_c_funct_bc1}. Writing $\lambda_n = (n + \rho(k))\eps$, we find that
\begin{equation*}
\frac{c^*_1(-\lambda_{n + 1}, k)}{c^*_1(-\lambda_{n + 1}, k + (0, 1))} = n, \quad 
\frac{c^*_s(-\lambda_n, k)}{c^*_s(-\lambda_n, k + (0, 1))} = n
\end{equation*}
for all $n \geq 1$, and similarly
\begin{align*}
\frac{c^*_1(-\lambda_{n + 1}, k)}{c^*_1(-\lambda_{n + 1}, k + (2, -1))} &= 4(n + k_2 + \tfrac{1}{2}), \qquad \; \frac{c^*_s(-\lambda_n, k)}{c^*_s(-\lambda_n, k + (2, -1))} = 4(n + k_2 - \tfrac{1}{2}), \\
\frac{\widetilde{c}_1(\lambda_{n + 1}, k + (0, -1))}{\widetilde{c}_1(\lambda_{n + 1}, k)} &= n + k_1 + 2k_2, \qquad \quad \;\; \frac{\widetilde{c}_s(\lambda_n, k + (0, -1))}{\widetilde{c}_s(\lambda_n, k)} = n + k_1 + 2k_2, \\
\frac{\widetilde{c}_1(\lambda_{n + 1}, k + (-2, 1))}{\widetilde{c}_1(\lambda_{n + 1}, k)} &= 4(n + k_1 + k_2 + \tfrac{1}{2}), \quad \frac{\widetilde{c}_s(\lambda_n, k + (-2, 1))}{\widetilde{c}_s(\lambda_n, k)} = 4(n + k_1 + k_2 - \tfrac{1}{2})
\end{align*}
for all $n \geq 0$, which match the corresponding shift factors from Lemma \ref{lem:nonsym_shift_factors} when taking the normalization of Definition \ref{def:nonsym_g_pm} into account.

The remaining cases are proven by induction. Suppose that the relation is valid for $\ell \in \bbZ_{\geq 0}B\backslash\{0\}$. Taking $\ell' \in B$, we compute
\begin{align*}
\eta(\lambda, \calG_+(\ell + \ell', k)) &= \eta(w^\lambda w_{\lambda^*}(\lambda^* - \rho(
\ell)\eps), \calG_+(k + \ell)) \eta(\lambda, \calG_+(\ell, k)) \\
&= \frac{\widetilde{c}_{w^\lambda w_{\lambda^*}}(\lambda^* - \rho(\ell)\eps + \rho(k + \ell)\eps, k + \ell + \ell')}{\widetilde{c}_{w^\lambda w_{\lambda^*}}(\lambda^* - \rho(\ell)\eps + \rho(k + \ell)\eps, k + \ell)} 
\frac{\widetilde{c}_{w^\lambda w_{\lambda^*}}(\lambda^* + \rho(k)\eps, k + \ell)}{\widetilde{c}_{w^\lambda w_{\lambda^*}}(\lambda^* + \rho(k)\eps, k)} \\
&= \frac{\widetilde{c}_{w^\lambda w_{\lambda^*}}(\lambda^* + \rho(k)\eps, k + \ell + \ell')}{\widetilde{c}_{w^\lambda w_{\lambda^*}}(\lambda^* + \rho(k)\eps, k)},
\end{align*}
so that the relation is also valid for $\ell + \ell'$. We can proceed similarly for $-\ell$ and $-\ell'$. This establishes the relation for all cases.
\end{proof}

\begin{proposition}\label{prop:ratio_non-symmetric_norms_bc1}
Let $\ell \in \bbZ_{\geq 0} B\backslash\{0\}$, $\lambda \in P$, $k \geq 1$, and $k - \ell \geq 0$. Then
\begin{equation*}
\frac{||E(\lambda, k)||_k^2}{||E(w^\lambda w_{\lambda^*}(\lambda^* + \rho(\ell)\eps), k - \ell)||_{k - \ell}^2} = \frac{\widetilde{c}_{w^\lambda w_{\lambda^*}}(\lambda^* + \rho(k)\eps, k - \ell) c^*_{w^\lambda w_{\lambda^*}}(-(\lambda^* + \rho(k)\eps), k)}{\widetilde{c}_{w^\lambda w_{\lambda^*}}(\lambda^* + \rho(k)\eps, k) c^*_{w^\lambda w_{\lambda^*}}(-(\lambda^* + \rho(k)\eps), k - \ell)}.
\end{equation*}
\end{proposition}

\begin{proof}
By Corollary \ref{cor:nonsym_shiftop_action_and_adjoint}, we have that
\begin{equation*}
\calG_+(\ell, k - \ell) E(w^\lambda w_{\lambda^*} (\lambda^* + \rho(\ell)\eps), k - \ell) = \eta(w^\lambda w_{\lambda^*} (\lambda^* + \rho(\ell)\eps), \calG_+(\ell, k - \ell)) E(\lambda, k),
\end{equation*}
and combining this with Corollary \ref{cor:nonsym_shiftop_action_and_adjoint}(ii) we find
\begin{align*}
&||E(\lambda, k)||_k^2 \\
&= \frac{1}{\eta(w^\lambda w_{\lambda^*} (\lambda^* + \rho(\ell)\eps), \calG_+(\ell, k - \ell))} (\calG_+(\ell, k - \ell) E(w^\lambda w_{\lambda^*} (\lambda^* + \rho(\ell)\eps), k - \ell), E(\lambda, k))_k \\
&= \frac{1}{\eta(w^\lambda w_{\lambda^*} (\lambda^* + \rho(\ell)\eps), \calG_+(\ell, k - \ell))} (E(w^\lambda w_{\lambda^*} (\lambda^* + \rho(\ell)\eps), k - \ell), \calG_-(-\ell, k)E(\lambda, k))_{k - \ell} \\
&= \frac{\eta(\lambda, \calG_-(-\ell, k))}{\eta(w^\lambda w_{\lambda^*} (\lambda^* + \rho(\ell)\eps), \calG_+(\ell, k - \ell))} ||E(w^\lambda w_{\lambda^*} (\lambda^* + \rho(\ell)\eps), k - \ell)||_{k - \ell}^2
\end{align*}
from which the statement follows after substituting the expressions from Proposition \ref{prop:shift_values_and_c_funct}.
\end{proof}

Finally, we recover the norm formula for the non-symmetric Jacobi polynomials.

\begin{theorem}\label{thm:rank_one_norm_bc1}
Let $k \in 2\bbZ \times \bbZ$ and suppose that $k \geq 1$. Then
\begin{equation*}
||E(\lambda, k)||_k^2 = \frac{c^*_{w^\lambda w_{\lambda^*}}(-(\lambda^* + \rho(k)\eps), k)}{\widetilde{c}_{w^\lambda w_{\lambda^*}}(\lambda^* + \rho(k)\eps, k)}, \quad \lambda \in P.
\end{equation*}
\end{theorem}

\begin{proof}
This follows from Proposition \ref{prop:ratio_non-symmetric_norms_bc1} by setting $k = \ell$, because $\widetilde{c}_{w^\lambda w_{\lambda^*}}(\lambda^* + \rho(k)\eps, 0) = 1$, and $c^*_{w^\lambda w_{\lambda^*}}(-(\lambda^* + \rho(k)\eps), 0) = 1$, and $||E(w^\lambda w_{\lambda^*}(\lambda^* + \rho(k)\eps), 0)||_0^2 = 1$.
\end{proof}


\subsection{Evaluation at the identity formula}
We view the polynomials $E(\lambda, k)$ as functions on $S^1$ by interpreting $z \in S^1$. To emphasize this we write
\begin{equation*}
E(\lambda, k; \cdot) : S^1 \to \bbC, \; z \mapsto E(\lambda, k; z).
\end{equation*}
Writing $e \in S^1$ for the identity element, the evaluation at the identity becomes the calculation of $E(\lambda, k; e)$. We establish this by invoking the methods of \cite{Opd89}, see also \cite[\S3.5]{HeSc94}, and adjusting it slightly to be compatible with $\calG_-(\ell, k)$.

The method hinges on the following property for symmetric lowering shift operators.

\begin{lemma}[{\cite[Cor.4.5]{Opd88a}}]\label{lem:sym_shift_rule_on_a}
For $\ell \in \bbZ_{\geq 0} B\backslash\{0\}$, we have
\begin{equation*}
G_-(-\ell, k)(f)(e) = G_-(-\ell, k)(1)(e) f(e), \quad f \in \bbC[S^1]^W.
\end{equation*}
\end{lemma}

There is a natural analog of this property for nonsymmetric lowering shift operators, which we prove using Lemma \ref{lem:sym_shift_rule_on_a}.

\begin{proposition}\label{prop:nonsym_shift_rule_on_a}
For $\ell \in \bbZ_{\geq 0} B\backslash\{0\}$, we have
\begin{equation*}
\calG_-(-\ell, k)(f)(e) = \calG_-(-\ell, k)(1)(e) f(e), \quad f \in \bbC[S^1].
\end{equation*}
\end{proposition}

\begin{proof}
Let $f \in \bbC[S^1]$ and note that $f_1, f_2 \in \bbC[S^1]^W$ by Lemma \ref{lem:steinberg_basis}. Using that Lemma \ref{lem:general_G_hat_explicit}, Lemma \ref{lem:sym_shift_rule_on_a}, and $f(e) = (f_1 + f_2 z)(e) = f_1(e) + f_2(e)$, we get
\begin{equation*}
\calG_-(-\ell, k)(f)(e) = G_-(-\ell, k_-)(f_1 + f_2)(e) = G_-(-\ell, k_-)(1)(e)f(e).
\end{equation*}
Taking $f = 1$ results in
\begin{equation*}
\calG_-(-\ell, k)(1)(e) = G_-(-\ell, k_-)(1)(e),
\end{equation*}
and substituting this in the previous equation gives the desired result.
\end{proof}

Now the proof of \cite[Thm.3.6.6]{HeSc94} can be modified to determine the evaluation at the identity for the non-symmetric Jacobi polynomials.

\begin{theorem}\label{thm:eval_at_e_nonsym_bc1}
Let $k \in 2\bbZ \times \bbZ$ and suppose that $k \geq 1$. Then
\begin{equation*}
E(\lambda, k; e) = \frac{\widetilde{c}_{w_0}(\rho(k)\eps, k)}{\widetilde{c}_{w^\lambda w_{\lambda^*}}(\lambda^* + \rho(k)\eps, k)}, \quad \lambda \in P.
\end{equation*}
\end{theorem}

\begin{proof}
Let $\ell \in \bbZ_{\geq 0}B \backslash \{ 0 \}$. By Proposition \ref{prop:shift_values_and_c_funct}, we have that
\begin{equation*}
\calG_-(-\ell, k) E(\lambda, k) = \frac{\widetilde{c}_{w^\lambda w_{\lambda^*}}(\lambda^* + \rho(k)\eps, k - \ell)}{\widetilde{c}_{w^\lambda w_{\lambda^*}}(\lambda^* + \rho(k)\eps, k)} E(w^\lambda w_{\lambda^*}(\lambda^* + \rho(\ell)\eps), k - \ell),
\end{equation*}
and specializing this to $\lambda = 0$ yields
\begin{equation*}
\calG_-(-\ell, k)(1) = \frac{\widetilde{c}_{w_0}(\rho(k)\eps, k - \ell)}{\widetilde{c}_{w_0}(\rho(k)\eps, k)} E(-\rho(\ell)\eps, k - \ell).
\end{equation*}
Combining this with Proposition \ref{prop:nonsym_shift_rule_on_a} yields
\begin{align*}
E(\lambda, k; e) &= \frac{\calG(-\ell, k)(E(\lambda, k))(e)}{\calG(-\ell, k)(1)(e)} \\
&= \frac{\widetilde{c}_{w^\lambda w_{\lambda^*}}(\lambda^* + \rho(k)\eps, k - \ell) \widetilde{c}_{w_0}(\rho(k)\eps, k)}{\widetilde{c}_{w^\lambda w_{\lambda^*}}(\lambda^* + \rho(k)\eps, k) \widetilde{c}_{w_0}(\rho(k)\eps, k - \ell)} \frac{E(w^\lambda w_{\lambda^*}(\lambda^* + \rho(\ell)\eps), k - \ell; e)}{E(-\rho(\ell)\eps, k - \ell; e)},
\end{align*}
as long as $E(-\rho(\ell)\eps, k - \ell; e) \neq 0$. By setting $\ell = k$, it reduces to
\begin{equation*}
E(\lambda, k; e) = \frac{\widetilde{c}_{w_0}(\rho(k)\eps, k)}{\widetilde{c}_{w^\lambda w_{\lambda^*}}(\lambda^* + \rho(k)\eps, k)},
\end{equation*}
as we have that $\widetilde{c}_{w^\lambda w_{\lambda^*}}(\lambda^* + \rho(k)\eps, 0) = 1$, $\widetilde{c}_{w_0}(\rho(k)\eps, 0) = 1$, $E(w^\lambda w_{\lambda^*}(\lambda^* + \rho(\ell)\eps), 0; e) = 1$, and $E(-\rho(\ell)\eps, 0; e) = 1$.
\end{proof}

\begin{corollary}\label{cor:eval_at_e_g_1}
Let $\ell \in \bbZ_{\geq 0}B \backslash \{ 0 \}$, $k \in 2\bbZ \times \bbZ$, and suppose that $k \geq 1$. Then
\begin{equation*}
\calG(-\ell, k)(1)(e) = \frac{\widetilde{c}_{w_0}(\rho(k - \ell)\eps, k - \ell)}{\widetilde{c}_{w_0}(\rho(k)\eps, k)}.
\end{equation*}
\end{corollary}

\begin{proofempty}
It follows immediately from the proof of Theorem \ref{thm:eval_at_e_nonsym_bc1} that
\pushQED{\qed}
\begin{equation*}
\calG(-\ell, k)(1)(e) = \frac{\widetilde{c}_{w_0}(\rho(k)\eps, k - \ell)}{\widetilde{c}_{w_0}(\rho(k)\eps, k)} \frac{\widetilde{c}_{w_0}(\rho(k - \ell)\eps, k - \ell)}{\widetilde{c}_{w_0}(\rho(\ell)\eps + \rho(k - \ell)\eps, k - \ell)} = \frac{\widetilde{c}_{w_0}(\rho(k - \ell)\eps, k - \ell)}{\widetilde{c}_{w_0}(\rho(k)\eps, k)}.\qedhere
\end{equation*}
\popQED
\end{proofempty}

\begin{remark}
Corollary \ref{cor:eval_at_e_g_1} is the natural analog of \cite[Cor.3.6.7]{HeSc94}, which states under the same assumptions that
\begin{equation*}
G_-(-\ell, k)(1)(e) = \frac{\widetilde{c}_1(\rho(k - \ell)\eps, k - \ell)}{\widetilde{c}_1(\rho(k)\eps, k)}.
\end{equation*}
\end{remark}


\section{Harish-Chandra homomorphism for non-symmetric Jacobi polynomials}\label{S:HC}


\subsection{Harish-Chandra homomorphism for symmetric Jacobi polynomials}
We recall from \cite{HeSc94} the required notions to state the main results on the Harish-Chandra homomorphism for the symmetric Jacobi polynomials (of type $BC_1$). For a complete account, also for general root systems, the reader should consult \cite[\S1]{HeSc94} or \cite{HeOp87}.

We cast some of the previous considerations into this formalism. Let $\fraka$ be $\bbR$ equipped with the standard inner product, and denote its complexification by $\frakh$. The root system $R = \{ \pm\eps, \pm2\eps \}$ and its weight lattice $P = \bbZ \eps$ are contained in $\fraka^*$. For each $\zeta \in \frakh$, we define a differential operator on $\bbC[P]$ by extending
\begin{equation}\label{eq:def_diff}
\partial_\zeta e^\lambda = \lambda(\zeta)e^\lambda, \quad \lambda \in P.
\end{equation}
Let $\xi$ be the unit vector of $\frakh$. Upon identifying $\bbC[P]$ with $\bbC[z^{\pm1}]$, we have that $\partial_\xi$ acts as $z\partial_z$, so we see that the modified Laplacian $ML_k$ from \eqref{eq:mod_lap} becomes
\begin{equation*}
ML_k = \partial_\xi^2 + k_1\frac{1 + e^{-\eps}}{1 - e^{-\eps}}\partial_\xi + 2k_2\frac{1 + e^{-2\eps}}{1 - e^{-2\eps}}\partial_\xi + \rho(k)^2.
\end{equation*}

We will consider an algebra of differential operators that contains $ML_k$. Let $S(\frakh)$ be the symmetric algebra associated with $\frakh$, and note that it may be identified with the algebra of polynomial functions on $\frakh^*$. Extending \eqref{eq:def_diff} yields an algebra of differential operators $S'(\frakh) = \{ \partial_p = p(\partial_\xi) \mid p \in S(\frakh) \}$. To allow for more general coefficients, which is needed for $ML_k$, we let $\frakR$ be the unital algebra generated by $(1 - e^{-\alpha})^{-1}$, with $\alpha \in R_+ = \{ \eps, 2\eps \}$, and define $\bbD_\frakR = \frakR \otimes S'(\frakh)$. The algebra $\bbD_\frakR$ acts on $\bbC_\Delta[P]$ by differential operators with coefficients in $\frakR$. Since $\bbZ_2$ acts on $\frakR$ and $S(\frakh)$, we can define an action of $\bbZ_2$ on $\bbD_\frakR$ by \begin{equation*}
w(f \otimes \partial_p) = w(f) \otimes \partial_{w(p)}, \quad f \otimes \partial_p \in \bbD_\frakR.
\end{equation*}
The subalgebra of $\bbZ_2$-invariant elements is denoted by $\bbD_\frakR^{\bbZ_2}$, which contains $ML_k$. By construction, the action of $\bbD_\frakR^{\bbZ_2}$ restricts to an action on $\bbC_\Delta[P]^{\bbZ_2}$.

There is also the algebra of differential-reflection operators $\bbD R_\frakR = \bbD_\frakR \otimes \bbC[\bbZ_2]$, with coefficients in $\frakR$, that contains the Cherednik operator $D_k$, because from \eqref{eq:cher_op} we see that 
\begin{equation*}
D_k = z\partial_z + k_1\frac{1}{1 - e^{-\eps}}(1 - s) + 2k_2\frac{1}{1 - e^{-2\eps}}(1 - s) - \rho(k).
\end{equation*}
The commutant of $1 \otimes \bbC[\bbZ_2]$ in $\bbD R_\frakR$ is denoted by $\bbD R_\frakR^{1 \otimes \bbC[\bbZ_2]}$, whose elements act on $\bbC_\Delta[P]^{\bbZ_2}$. There is a well-defined algebra homomorphism
\begin{equation}\label{eq:def_beta}
\beta : \bbD R_\frakR^{1 \otimes \bbC[\bbZ_2]} \to \bbD_\frakR^{\bbZ_2}, \quad D \otimes w \mapsto D,
\end{equation}
see \cite[Lem.1.2.2]{HeSc94}, with the property for each $\calD \in \bbD R_\frakR^{1 \otimes \bbC[\bbZ_2]}$
\begin{equation*}
\beta(\calD)f = \calD f, \quad f \in \bbC_\Delta[P]^{\bbZ_2}.
\end{equation*}
In particular, the element $D_k^2 \in \bbD R_\frakR^{1 \otimes \bbC[\bbZ_2]}$ is mapped onto $ML_k$, cf. \S\ref{S:shiftops_for_sym_jac}.

\begin{definition}\label{def:scalar_const_term_and_hc}
The constant term map $\gamma'(k)$ and the Harish-Chandra homomorphism $\gamma(k)$ are the unique algebra homomorphisms from $\bbD_\frakR$ to $S(\frakh)$ satisfying
\begin{align*}
\gamma'(k) \left( \frac{1}{1 - e^{-\alpha}} \right) &= 1, \quad \gamma(k) \left( \frac{1}{1 - e^{-\alpha}} \right) = 1, \quad \alpha \in R_+, \\
\quad \gamma'(k)(\partial_\xi) &= \xi, \qquad \qquad \; \thinspace \gamma(k)(\partial_\xi) = \xi - \rho(k)\eps.
\end{align*}
\end{definition}

Let us summarize the main properties of the constant term map and the Harish-Chandra homomorphism that we seek to replicate with their analogs for the non-symmetric Jacobi polynomials.

\begin{theorem}[{\cite[Thm.1.3.12]{HeSc94}, \cite[Thm.2.12(2)]{Opd95}}]\text{ }\label{thm:sym_hc}
\begin{enumerate}[(i)]
    \item If $D \in \bbD_\frakR$ has the polynomials $p(\lambda, k)$ as eigenfunctions, then
    \begin{equation*}
    D p(\lambda, k) = \gamma'(k)(D)(\lambda)p(\lambda, k), \quad \lambda \in P^+.
    \end{equation*}
    \item Let $\mathbb{D}(k)$ be the commutant of $ML_k$ in $\bbD_\frakR^{\bbZ_2}$. Then
    \begin{equation*}
    \gamma(k) : \mathbb{D}(k) \to S(\frakh)^{\bbZ_2}
    \end{equation*}
    is a well-defined algebra isomorphism.
    \item For $D \in \bbD(k)$ and $p = \gamma(k)(D)$, we have $D = \beta(p(D_k))$, and
    \begin{equation*}
    D p(\lambda, k) = p(\widetilde{\lambda})p(\lambda, k), \quad \lambda \in P^+.
    \end{equation*}
    \end{enumerate}
\end{theorem}


\subsection{Non-symmetric Jacobi polynomials as $\bbC^2$-valued Laurent polynomials}
At the heart of the previous discussion lies the fact that the action of the reflection can be absorbed by the invariance of the symmetric Jacobi polynomials. This allows us to restrict the action of differential-reflection operators to just differential operators, for which a constant term map is available. If we want to replicate this for the non-symmetric Jacobi polynomials, then we should realize the non-symmetric Jacobi polynomials as invariant vector-valued Laurent polynomials.

We recall and extend some notions from \cite[\S2.1]{vHvP1}. The map from \cite[(2)]{vHvP1} can be extended to a $\bbC_\Delta[P]$-module isomorphism
\begin{equation*}
\Gamma : \bbC_\Delta[P] \to (\bbC_\Delta[P] \otimes \bbC^2)^{\bbZ_2}, \quad p \mapsto (p, s(p))^T.
\end{equation*}
Note that $\Gamma(\bbC[P]) = (\bbC[P] \otimes \bbC^2)^{\bbZ_2}$, and recall that $\bbC[P] \otimes \bbC^2$ is equipped with a sesquilinear pairing such that $\Gamma|_{\bbC[P]} : \bbC[P] \to (\bbC[P] \otimes \bbC^2)^{\bbZ_2}$ respects the sesquilinear pairings. Explicitly,
\begin{equation*}
(P, Q)_k = \frac{1}{2}\int_{S^1}\left(\overline{P_1(z)}Q_1(z) + \overline{P_2(z)} Q_2(z)\right)\delta_k(z) \frac{dz}{iz}.
\end{equation*}
Let $P(\lambda, k) = \Gamma(E(\lambda, k))$, with $\lambda \in P$, be the vector-valued Jacobi polynomials. Then $(P(\lambda, k) \mid \lambda \in P)$ is an orthogonal basis of $(\bbC[P] \otimes \bbC^2)^{\bbZ_2}$.

To formulate the push-forward $\Gamma_*$ of $\Gamma$ properly, we introduce the algebra $\bbD R_\frakR \otimes \bbM_2$, whose product is given by
\begin{equation*}
\calD \otimes A \cdot \calE \otimes B = \calD \calE \otimes AB.
\end{equation*}
Let $\Delta : \bbC[\bbZ_2] \to \bbD R_\frakR \otimes \bbM_2, \; w \mapsto 1 \otimes w \otimes \pi(w)$. We denote the commutant of $\Delta(\bbC[\bbZ_2])$ in $\bbD R_\frakR \otimes \bbM_2$ by $(\bbD R_\frakR \otimes \bbM_2)^{\Delta(\bbC[\bbZ_2])}$, whose elements act on $\left(\bbC_\Delta[P]\otimes\bbC^2\right)^{\bbZ_2}$.

\begin{lemma}[{\cite[(3)]{MvP-VVOP}}]
The map
\begin{equation*}
\Gamma_* : \bbD R_\frakR \to (\bbD R_\frakR \otimes \bbM_2)^{\Delta(\bbC[\bbZ_2])}, \quad \calD \mapsto \diag(\calD, s \calD s).
\end{equation*}
is a well-defined algebra homomorphism. Moreover, for each $\calD \in \bbD R_\frakR$
\begin{equation*}
\Gamma_*(\calD) \Gamma(f) = \Gamma(\calD f), \quad f \in \bbC_\Delta[P].
\end{equation*}
\end{lemma}

Similarly, we define the algebra $\bbD_\frakR \otimes \bbM_2$. Let $\pi : \bbZ_2 \to \bbM_2 = \End(\bbC^2)$ be the representation of $\bbZ_2$ on $\bbC^2$ such that
\begin{equation}\label{eq:pi_s}
\pi(s) = \begin{pmatrix} 0 & 1 \\ 1 & 0 \end{pmatrix}.
\end{equation}
The action of $\bbZ_2$ on $\bbD_\frakR \otimes \bbM_2$ is given by
\begin{equation*}
w(D \otimes A) = w(D) \otimes \pi(w) A \pi(w^{-1}), \quad D \otimes A \in \bbD_\frakR \otimes \bbM_2.
\end{equation*}
The subalgebra of $\bbZ_2$-invariant elements is denoted by $\bbD_\frakR^{\bbZ_2}$, whose elements act on $\bbC_\Delta[P]^{\bbZ_2}$.

Finally, we define the matrix-valued analog $\beta_{\ns}$ of the map $\beta$ from \eqref{eq:def_beta}.

\begin{lemma}
The map
\begin{equation*}
\beta_\ns : (\bbD R_\frakR \otimes \bbM_2)^{\Delta(\bbC[\bbZ_2])} \to (\bbD_\frakR \otimes \bbM_2)^{\bbZ_2}, \; D \otimes w \otimes A \mapsto D \otimes A \pi(w^{-1})
\end{equation*}
is a well-defined algebra homomorphism. Moreover, for each $\calD \in (\bbD R_\frakR \otimes \bbM_2)^{\Delta(\bbC[\bbZ_2])}$
\begin{equation*}
\beta_\ns(\calD) \calF = \calD \calF, \quad \calF \in (\bbC_\Delta[P] \otimes \bbC^2)^{\bbZ_2}.
\end{equation*}
\end{lemma}

\begin{proof}
Consider the map $\Delta' : \bbC[\bbZ_2] \to \bbD_\frakR \otimes \bbM_2 \otimes \bbC[\bbZ_2], \; w \mapsto 1 \otimes 1 \otimes w$ and note that
\begin{equation*}
\chi : (\bbD R_\frakR \otimes \bbM_2)^{\Delta(\bbC[\bbZ_2])} \to (\bbD_\frakR \otimes \bbM_2 \otimes \bbC[\bbZ_2])^{\Delta'(\bbC[\bbZ_2])}, \quad D \otimes w \otimes A \mapsto D \otimes A \pi(w^{-1}) \otimes w
\end{equation*}
is a well-defined algebra homomorphism. We recall from \cite[Lem.5.3]{MvP-VVOP} that
\begin{equation*}
\beta'_\ns : (\bbD_\frakR \otimes \bbM_2 \otimes \bbC[\bbZ_2])^{\Delta'(\bbC[\bbZ_2])} \to (\bbD_\frakR \otimes \bbM_2)^{\bbZ_2}, \quad D \otimes A \otimes w \mapsto D \otimes A
\end{equation*}
is also a well-defined algebra homomorphism. It follows that their composition
\begin{equation*}
\beta_\ns = \beta'_\ns \circ \chi : (\bbD R_\frakR \otimes \bbM_2)^{\Delta(\bbC[\bbZ_2])} \to (\bbD_\frakR \otimes \bbM_2)^{\bbZ_2}
\end{equation*}
is a well-defined algebra homomorphism as well.

The second statement follows from the observation that the actions of $\diag(s, s)$ and $\pi(s)$ agree on $(\bbC_\Delta[P] \otimes \bbC^2)^{\bbZ_2}$.
\end{proof}

\begin{example}\label{ex:ct_cher}
Applying the previous two maps to the Cherednik operator $D_k$, we recover
\begin{equation*}
\widetilde{D}_k = \beta_\ns(\Gamma_*(D_k)) =
\begin{pmatrix}
\partial_\xi + \frac{k_1}{1 - e^{-\eps}} + \frac{2k_2}{1 - e^{-2\eps}} - \rho(k) & -\frac{k_1}{1 - e^{-\eps}} - \frac{2k_2}{1 - e^{-2\eps}} \\
-\frac{k_1}{1 - e^{\eps}} - \frac{2k_2}{1 - e^{2\eps}} & -\partial_\xi - \frac{k_1}{1 - e^{-\eps}} - \frac{2k_2}{1 - e^{-2\eps}} + \rho(k)
\end{pmatrix},
\end{equation*}
see \cite[(4)]{vHvP1}.
\end{example}


\subsection{$\bbM_2$-valued Harish-Chandra homomorphism}
We define the matrix-valued analogs of the constant term map and the Harish-Chandra homomorphism, that act on elements of $\bbD_\frakR \otimes \bbM_2$, and study their properties. In particular, an analog is established for each of the three statements from Theorem \ref{thm:sym_hc}.

\begin{definition}
The matrix-valued constant term map $\gamma'_\ns(k)$ and the matrix-valued Harish-Chandra homomorphism $\gamma_\ns(k)$ are the unique algebra homomorphisms from $\bbD_\frakR \otimes \bbM_2$ to $S(\frakh) \otimes \bbM_2$ satisfying
\begin{equation*}
\gamma_\ns'(k)(D \otimes A) = \gamma'(k)(P) \otimes A, \quad \gamma_\ns(k)(D \otimes A) = \gamma(k)(D) \otimes A
\end{equation*}
for all $D \otimes A \in \bbD_\frakR \otimes \bbM_2$.
\end{definition}

\begin{example}
The constant term of $\widetilde{D}_k$ is given by
\begin{equation*}
\gamma'_\ns(\widetilde{D}_k) =
\begin{pmatrix}
\partial_\xi + \rho(k) & -2\rho(k) \\
0 & -\partial_\xi - \rho(k)
\end{pmatrix}.
\end{equation*}
\end{example}

\begin{proposition}\label{prop:mv_hc_inj}
Let $\widetilde{\bbD}_\ns(k)$ be the commutant of $\widetilde{D}_k$ in $\bbD_\frakR \otimes \bbM_2$.
\begin{enumerate}[(i)]
\item The restriction
\begin{equation*}
\gamma_\ns'(k) : \widetilde{\bbD}_\ns(k) \to S(\frakh) \otimes \bbM_2
\end{equation*}
is an injective algebra homomorphism.
\item The image $\gamma_\ns'(k)(\widetilde{\bbD}_\ns(k))$ is contained in the upper triangular matrices.
\item For $D \in \widetilde{\bbD}_\ns(k)$, the entry $\gamma_\ns'(k)(D)_{12}$ is determined by $\gamma_\ns'(k)(D)_{11}$ and $\gamma_\ns'(k)(D)_{22}$.
\end{enumerate}
The properties (i)--(iii) also hold for the restriction of $\gamma_\ns(k)$ to $\widetilde{\bbD}_\ns(k)$.
\end{proposition}

\begin{proof}
Let $D \in \widetilde{\bbD}_\ns(k)$, i.e., $D \in \bbD_\frakR \otimes \bbM_2$ commutes with $\widetilde{D}_k$. We expand the matrix entries of $D$ and $\widetilde{D}_k$ as follows
\begin{equation*}
D_{ij} = \sum_{n \geq 0} e^{-n\eps} \partial_{P_n^{ij}}, \quad (\widetilde{D}_k)_{ij} = \sum_{m \geq 0} e^{-m\eps} \partial_{Q_m^{ij}}, \quad P_n^{ij}, Q_m^{ij} \in S(\frakh).
\end{equation*}
Using the rule $\partial_{p(\lambda)} \circ e^{-n\eps} = e^{-n\eps} \partial_{p(\lambda - n\eps)}$ and that $D$ commutes with $\widetilde{D}_k$, we find that
\begin{equation*}
\sum_{n \geq 0} e^{-n\eps} \sum_{0 \leq m \leq n, \; \ell = 1,2} \left(\partial_{P_{n - m}^{i\ell}(\lambda - m\eps) Q_m^{\ell j}(\lambda) - Q_m^{i\ell}(\lambda - (n - m)\eps)P_{n - m}^{\ell j}(\lambda)}\right) = 0.
\end{equation*}
Restricting to the coefficient of $e^{-n\eps}$ yields
\begin{equation*}
\sum_{\ell = 1}^2 \left(P_n^{i\ell}(\lambda)Q_0^{\ell j}(\lambda) - Q_0^{i\ell}(\lambda - n\eps)P_n^{\ell j}(\lambda)\right) + \Sigma_n^{ij}(\lambda) = 0,
\end{equation*}
with
\begin{equation*}
\Sigma_n^{ij}(\lambda) = \sum_{0 < m \leq n, \; \ell = 1,2} \left(P_{n - m}^{i\ell}(\lambda - m\eps) Q_m^{\ell j}(\lambda) - Q_m^{i\ell}(\lambda - (n - m)\eps)P_{n - m}^{\ell j}(\lambda)\right).
\end{equation*}
Using Example \ref{ex:ct_cher}, we obtain four explicit expressions
\begin{equation}\label{eq:zar_eqs}
\begin{aligned}
\xi(-n\eps) P_n^{11}(\lambda) &= \Sigma_n^{11}(\lambda) + \rho(k) P^{21}_n(\lambda), \\
(\xi(2\lambda - n\eps) + 2\rho(k)) P_n^{12}(\lambda) &= \Sigma_n^{12}(\lambda) - \rho(k) (P^{11}_n(\lambda) - P^{22}_n(\lambda)), \\
(\xi(2\lambda - n\eps) + 2\rho(k)) P_n^{21}(\lambda) &= \Sigma_n^{21}(\lambda), \\
-\xi(-n\eps) P_n^{22}(\lambda) &= \Sigma_n^{22}(\lambda) - \rho(k) P^{21}_n(\lambda).
\end{aligned}
\end{equation}
The first and fourth equations from \eqref{eq:zar_eqs} determine the polynomials $P^{11}_n$ and $P^{22}_n$ as long as $n \geq 1$, and the second and third equations from \eqref{eq:zar_eqs} determine the polynomials $P^{12}_n$ and $P^{21}_n$ on the Zariski dense set
\begin{equation*}
\{ \lambda \in \frakh^* \mid \xi(2\lambda - n\eps) + 2\rho(k) \neq 0 \},
\end{equation*}
which is sufficient to recover the polynomials completely. It follows that $(P_n^{ij})_{ij}$ is determined by $\{ (P_m^{ij})_{ij} \}_{0 \leq m < n}$. By induction, we obtain that $D$ can be recovered from its constant term, which proves (i).

For $n = 0$, the second and third equations from \eqref{eq:zar_eqs} still determine $P_0^{12}$ and $P_0^{21}$. Note that $\Sigma_0^{ij} = 0$, so that $P_0^{21} = 0$ and $P_0^{12}$ is determined by $P_0^{11}$ and $P_0^{22}$. This settles (ii) and (iii). As for the last statement, we recall from \cite[(1.2.9)]{HeSc94} that, formally, $\gamma(k)(D) = e^{\rho(k)\eps} \circ \gamma'(k)(D) \circ e^{-\rho(k)\eps}$ for $D \in \bbD_\frakR$. This allows us to transfer the statements.
\end{proof}

We prove the analog of Theorem \ref{thm:sym_hc}(i) for the matrix-valued constant term map.

\begin{corollary}\label{cor:mv_ct_eigenvalue_prop}
If $D \in \bbD_\frakR \otimes \bbM_2$ has the polynomials $P(\lambda, k)$ as eigenfunctions, then
\begin{align*}
D P(n + 1, k) &= \ct(D)_{11}((n + 1)\eps) P(n + 1, k), \\
D P(-n, k) &= \ct(D)_{22}(n\eps) P(-n, k)
\end{align*}
for all $n \in \bbZ_{\geq 0}$.
\end{corollary}

\begin{proof}
Since the polynomials $P(\lambda, k)$ are a basis of simultaneous eigenfunctions of the operators $D$ commutes with $\widetilde{D}_k$, we find that $D$ commutes with $\widetilde{D}_k$, i.e., $D \in \widetilde{\bbD}_\ns(k)$.

For $n \in \bbZ_{\geq 0}$, we have
\begin{align*}
P(n + 1, k) &= 
\begin{pmatrix}
e^{(n + 1)\eps} \\
e^{-(n + 1)\eps}
\end{pmatrix}
+ \text{lower order terms}, \\
P(-n, k) &= 
\begin{pmatrix}
e^{-n\eps} + d_n(k) e^{n\eps} \\
d_n(k) e^{-n\eps} + e^{n\eps}
\end{pmatrix}
+ \text{lower order terms}, \quad d_n(k) \in \bbC.
\end{align*}
If we expand the entries of $D$ as in the proof of Proposition \ref{prop:mv_hc_inj}, then the eigenvalue can be read off from the $e^{(n + 1)\eps}$ and $e^{n\eps}$ terms. In the case of $P(n + 1, k)$, we consider the first component and find that the only contribution comes from
\begin{equation*}
\partial_{\ct(D)_{11}}(e^{(n + 1)\eps}) = \ct(D)_{11}((n + 1)\eps)e^{(n + 1)\eps},
\end{equation*}
which gives the stated eigenvalue. For $P(-n, k)$ we need to look at the second component and use that $\ct(D)_{21} = 0$, see Proposition \ref{prop:mv_hc_inj}.
\end{proof}

Let $\diag : S(\frakh) \otimes \bbM_2 \to S(\frakh) \otimes \bbC^2, P \mapsto (P_{11}, P_{22})$ be the diagonal map.

\begin{corollary}[of Proposition \ref{prop:mv_hc_inj}]\label{cor:comp_diag_hc}
The composition
\begin{equation*}
\diag \circ \gamma_\ns(k) : \widetilde{\bbD}_\ns(k) \to S(\frakh) \otimes \bbC^2
\end{equation*}
is a well-defined injective algebra homomorphism.
\end{corollary}

Note that it is not immediately clear whether the image of the map from Corollary \ref{cor:comp_diag_hc} is invariant for some action of $\bbZ_2$. We calculate the images of the operators of the form $\beta_\ns(\Gamma_*(p(D_k)))$, with $p \in S(\frakh)$.

\begin{lemma}\label{lem:CT_of_p}
For $p \in S(\frakh)$, we have
\begin{equation*}
(\gamma_\ns(k) \circ \beta_\ns \circ \Gamma_*)(p(D_k)) = \begin{pmatrix} p & \rho(k)\frac{s(p) - p}{\xi} \\ 0 & s(p) \end{pmatrix}.
\end{equation*}
\end{lemma}

\begin{proof}
Both expressions define algebra homomorphisms from $S(\frakh)$ to $S(\frakh) \otimes \bbM_2$. Since $S(\frakh)$ is generated by $1$ and $\xi$, it is sufficient to check whether the expressions coincide for $1$ and $\xi$, the latter follows directly from Example \ref{ex:ct_cher}.
\end{proof}

Because of Lemma \ref{lem:CT_of_p}, we define the action of $\bbZ_2$ on $S(\frakh) \otimes \bbC^2$ by
\begin{equation*}
s(P) = (s(P_2), s(P_1))^T, \quad P \in S(\frakh) \otimes \bbC^2.
\end{equation*}
Before we can show this invariance, albeit on a subalgebra of $\widetilde{\bbD}_\ns(k)$, we prove the following.

\begin{lemma}\label{lem:scalar_hc_bc1}
Any element $\calD \in \bbD R_\frakR$ that commutes with $D_k$ is of the form
\begin{equation*}
\calD = p(D_k), \quad p \in S(\frakh).
\end{equation*}
\end{lemma}

\begin{proof}
We seek to apply the proof of Proposition \ref{prop:pre_nonsym_shift_str_thm}. For this we rewrite it in terms of the coordinate $z$, we have $\partial_\xi = z\partial_z$ and $\calD \in \bbC(\tfrac{1}{1 - z^{-1}}) \otimes \bbC[z\partial_z] \otimes \bbC[\bbZ_2]$ satisfying $[D_k, \calD] = 0$. Following the proof of Proposition \ref{prop:pre_nonsym_shift_str_thm}, we find
\begin{equation*}
\calD = a_N(k) (z\partial_z)^N + \sum_{i \leq N - 1, \; j = 0, 1} c_{ij}(k) (z\partial_z)^i s^j, \quad c_{ij}(k) \in \bbC(\tfrac{1}{1 - z^{-1}}),
\end{equation*}
and $a_N(k) \in \bbC(\tfrac{1}{1 - z^{-1}})$ such that $z\partial_z(a_N(k)) = 0$. This implies that $a_N(k)$ is constant, i.e., $a_N(k) \in \bbC$. In particular, we have that each element of $\bbD R_{\frakR}$ admits an order. The statement now follows by induction on the order of $\calD$, cf. the proof of Corollary \ref{cor:pre_nonsym_shift_str_thm}.
\end{proof}

Let $\pi_1 : S(\frakh) \otimes \bbC^2 \to S(\frakh)$ be the projection onto the first component. We note that $\pi_1 : (S(\frakh) \otimes \bbC^2)^{\bbZ_2} \to S(\frakh)$ is an algebra isomorphism. We are in the position to prove the analog of Theorem \ref{thm:sym_hc}(ii).

\begin{theorem}\label{thm:HCinv_BC1_and_A1}
Let $\bbD_\ns(k)$ be the commutant of $\widetilde{D}_k$ in $(\bbD_\frakR \otimes \bbM_2)^{\bbZ_2}$. Then
\begin{equation*}
\diag \circ \gamma_\ns(k): \bbD_\ns(k) \to (S(\frakh) \otimes \bbC^2)^{\bbZ_2}
\end{equation*}
is a well-defined algebra isomorphism. Moreover, the composition
\begin{equation*}
\gamma_{\mathrm{HC}}(k) := \pi_1 \circ (\diag \circ \gamma_\ns(k)): \bbD_\ns(k) \to S(\frakh)
\end{equation*}
is an algebra isomorphism that we call the (non-symmetric) Harish-Chandra isomorphism.
\end{theorem}

\begin{proof}
Let $D \in \bbD_\ns(k)$, i.e., $D$ is $\bbZ_2$-invariant and commutes with $\widetilde{D}_k$. Let $E_{ij}$ denote the matrix with $(i, j)$-entry $1$ and all other entries $0$. Note that $\pi(s)E_{11}\pi(s^{-1}) = E_{22}$ and $\pi(s)E_{12}\pi(s^{-1}) = E_{21}$, so the $s$-invariance of $D$ implies that
\begin{equation*}
s(D_{11}) = D_{22}, \quad s(D_{12}) = D_{21}.
\end{equation*}
From Example \ref{ex:ct_cher}, we see that $D_k = (\widetilde{D}_k)_{11} + (\widetilde{D}_k)_{12}s$. Similarly, we define $\calD = D_{11} + D_{12}s$. Note that $D_k$ and $\calD$ are both elements of $\bbD R_{\frakR}$. We compute
\begin{align*}
[D_k, \calD] &= ((\widetilde{D}_k)_{11}D_{11} - D_{11}(\widetilde{D}_k)_{11} + (\widetilde{D}_k)_{12}D_{21} - D_{12}(\widetilde{D}_k)_{21}) \\
&+ ((\widetilde{D}_k)_{11}D_{12} - D_{12}(\widetilde{D}_k)_{22} + (\widetilde{D}_k)_{12}D_{22} - D_{11}(\widetilde{D}_k)_{12})s = 0,
\end{align*}
which follows from comparing the $(1,1)$- and $(1,2)$-entries of the relation $[\widetilde{D}_k, D] = 0$.

By Lemma \ref{lem:scalar_hc_bc1}, there exists a unique polynomial $p \in S(\frakh)$ such that $\calD = p(D_k)$.
\begin{align*}
\beta_\ns(\Gamma_*(\calD)) &= \beta_\ns(\diag(D_{11}, s(D_{11}))) + \beta_\ns(\diag(D_{12}s, s(D_{12})s)) \\
&= \begin{pmatrix} D_{11} & D_{12} \\ s(D_{12}) & s(D_{11}) \end{pmatrix} = D.
\end{align*}
At the same time, Lemma \ref{lem:CT_of_p} yields
\begin{equation*}
\diag(\gamma_\ns(k)(\beta_\ns(\Gamma_*(\calD))) = \diag(p, s(p)) \in (S(\frakh) \otimes \bbC^2)^{\bbZ_2}.
\end{equation*}

Since any element of $(S(\frakh) \otimes \bbC^2)^{\bbZ_2}$ is of the form $\diag(p, s(p))$, with $p \in S(\frakh)$, and because $\beta_\ns(\Gamma_*(p(D_k)) \in \bbD_\ns(k)$, the surjectivity of $\diag \circ \gamma_\ns$ also follows from this.
\end{proof}

From the proof of Theorem \ref{thm:HCinv_BC1_and_A1}, we also obtain the analog of Theorem \ref{thm:sym_hc}(iii).

\begin{corollary}\label{cor:H-C-homEigenvalue}
For $D \in \bbD_\ns(k)$ and $p = \gamma_\ns(k)(D)_{11}$, we have $D = \beta_\ns(\Gamma_*(p(D_k)))$, and
\begin{equation*}
D \calP(\lambda, k) = p(\widetilde{\lambda})\calP(\lambda, k), \quad \lambda \in P.
\end{equation*}
\end{corollary}

\begin{proof}
This follows from the proof of Theorem \ref{thm:HCinv_BC1_and_A1}, together with \eqref{eq:BC1_eigenvalue}.
\end{proof}


\section{Relation to spherical functions}\label{S:SF}


\subsection{Eigenvalue map}\label{S:Eigenvalue_map}
The dependency of the eigenvalue $p(\widetilde{\lambda})$ from Corollary \ref{cor:H-C-homEigenvalue} on the spectral parameter $\lambda \in P = \bbZ\eps$ is not polynomial. To remedy this in a formal way, we collect the spectral parameters $-N\eps$ and $(N + 1)\eps$ together into the pair $(-N, N + 1)$, where $N \in \bbZ_{\geq 0}$. We combine the  $\bbC^2$-valued Laurent polynomials for these spectral parameters into $\bbM_2$-valued Laurent polynomials.

\begin{definition}
The $\bbM_2$-valued Laurent polynomial $M(N, k)$ with spectral parameter $N \in \bbZ_{\geq 0}$ is defined as the matrix whose columns are $P(-N, k)$ and $P(N + 1, k)$.
\end{definition}

The polynomials $M(N, k)$ are invariant for the action of $\bbZ_2$ on $\bbC[P] \otimes \bbM_2$, given by
\begin{equation*}
(s\cdot M)(z) = \pi(s)M(s^{-1}(z)).
\end{equation*}
The algebra $\bbD_\ns(k)$ acts on $\left(\bbC[P] \otimes \bbM_2\right)^{\bbZ_2}$ by its natural action on the columns of the $\bbM_2$-valued functions.

\begin{definition}
The eigenvalue map $\gamma_{\mathrm{EV}} : \bbD_\ns(k) \to S(\frakh) \otimes \bbM_2$ is defined by
\begin{equation*}
\gamma_{\mathrm{EV}}(D)(N\eps) =
\begin{pmatrix}
\gamma_{\mathrm{HC}}(k)(D)(-N\eps - \rho(k)\eps) & 0 \\
0 & \gamma_{\mathrm{HC}}(k)(D)((N + 1)\eps + \rho(k)\eps)
\end{pmatrix},
\quad N \in \bbZ.
\end{equation*}
\end{definition}
Note that $\gamma_{\mathrm{EV}}$ is an algebra homomorphism. The polynomials $M(N,k)$ are not eigenfunctions for the action of $D \in \bbD_\ns(k)$, however, we do have
\begin{equation}\label{eq:D_acting_on_M}
DM(N,k) = M(N,k) \gamma_{\mathrm{EV}}(D)(N\eps), \quad N \in \bbZ_{\geq 0}.
\end{equation}
In this sense we say that the polynomials $M(N,k)$ are eigenfunctions of $D$ with diagonal eigenvalue $\gamma_{\mathrm{EV}}(D)(N\eps)$. Moreover, the dependency of the eigenvalues on the spectral parameter $N$ is polynomial. Note that for a polynomial $p$ in one variable, which we can also view as an element of $S(\lah)$, we have
\begin{equation}\label{eq:eigenvalue_map_pD}
\gamma_{\mathrm{EV}}(\beta_\ns(\Gamma_*(p(D_k)))(N\eps) =
\begin{pmatrix}
p(-N\eps - \rho(k)\eps) & 0 \\
0 & p((N + 1)\eps + \rho(k)\eps)
\end{pmatrix}, \quad N \in \bbZ_{\geq 0}.
\end{equation}
For $s\in\bbZ_2$ we define $s_*(\lambda)=s(\lambda+\rho(k)\eps)-\rho(k)\eps$. We observe that the function
\begin{equation*}
\lambda + \tfrac{1}{2}\eps \mapsto\gamma_{\mathrm{EV}}(\beta_\ns(\Gamma_*(p(D_k)))(\lambda),
\end{equation*}
viewed as an element of $S(\lah) \otimes \bbM_2$, is invariant for the action of $\bbZ_2$ on $S(\lah) \otimes \bbM_2$, given by
\begin{equation}\label{eqn:s_*_action}
(s \cdot (p \otimes A))(\lambda) = p(s_*^{-1}\lambda) \otimes \left(\pi(s) A \pi(s)^{-1}\right).
\end{equation}


\subsection{Spherical functions}
Recall from \cite[\S3]{vHvP1} that the spherical functions $\Psi_\pm(\ell)$, $\ell \in \bbZ_{\geq 0}$, for the compact symmetric pair $(\Spin(2m+2),\Spin(2m+1))$ of type $\tau$ can be identified with the non-symmetric Jacobi polynomials. The spherical functions are simultaneous eigenfunctions of the algebra $\bbC[R_m]$, where $R_m$ is, up to a scalar multiple, the radial part of the Dirac operator, see \cite[Lem.3.1]{vHvP1}. In particular, we have
\begin{equation*}
R_m\Psi_\pm(\ell) = \pm(2\ell + 2m+ 1)\Psi_\pm(\ell), \quad \ell \in \bbZ_{\geq 0}.
\end{equation*}

To discuss the dependence of the eigenvalues on the spectral parameters $(\sigma^\pm, \ell)$ in the representation-theoretic context, we combine the spherical functions, that are $\bbC^2$-valued Laurent polynomials, into $\bbM_2$-valued Laurent polynomials. 

\begin{definition}
The full spherical function $\Psi(\ell)$ with spectral parameter $\ell$ is defined as the $\bbM_2$-valued Laurent polynomial whose columns are $\Psi_-(\ell)$ and $\Psi_+(\ell)$.
\end{definition}

Similarly to \eqref{eq:D_acting_on_M} and \eqref{eq:eigenvalue_map_pD}, we have for a polynomial $p$ in one variable that
\begin{equation*}
p(R_m) \Psi(\ell) = \Psi(\ell)
\begin{pmatrix}
p(-(2\ell + 2m + 1)\eps) & 0 \\
0 & p((2\ell + 2m + 1)\eps)
\end{pmatrix},
\quad \ell \in \bbZ_{\geq 0},
\end{equation*}
and that the function
\begin{equation*}
\ell + \tfrac{1}{2} \mapsto
\begin{pmatrix}
p(-(2\ell + 2m + 1)\eps) & 0 \\
0 & p((2\ell + 2m + 1)\eps)
\end{pmatrix}
\end{equation*}
is invariant for the action \eqref{eqn:s_*_action} of $\bbZ_2$. We observe that the full spherical functions $\Psi(\ell)$ have eigenvalues that are polynomials in the spectral parameter $\ell$, and that they have $\bbZ_2$-invariance similar to $\gamma_{\mathrm{EV}}(D)$, with $D \in \bbD_\ns(k)$.

We will give an alternative explanation for these properties based on the representation-theoretic context as discussed in \cite[Ch.3]{MR0929683}. First, we will have to introduce some notation. Although the following results hold for $m \geq 1$, we exclude the case $m = 1$. The reason for this is that $\Spin(4)$ is not simple, which would require us to argue separately for this case throughout the next section.

Let $m > 1$ and let $\lag = \laso(2m+2,\bbC)$, $\lak = \laso(2m+1,\bbC)$, and $\lam=\laso(2m,\bbC)$ be the complexified Lie algebras of the compact spin groups
\begin{equation*}
U_m = \Spin(2m+2), \quad K_m = \Spin(2m+1), \quad M_m = \Spin(2m),
\end{equation*}
respectively. Note that we drop the subscript $m$ from the notation for the Lie algebras. We follow \cite{MR0240238} for our choices of maximal tori, roots, weights, and positivity. The maximal tori of $\lag$ and $\lam$ are denoted by $\lat$ and $\lat_\lam$, respectively. The fundamental weights of $U_m$, $K_m$, and $M_m$, or equivalently of $\lag$, $\lak$, and $\lam$, are denoted by
\begin{equation*}
\varpi_1,\ldots,\varpi_{m+1}, \quad\omega_1,\ldots,\omega_m, \quad\eta_1,\ldots, \eta_m,
\end{equation*}
respectively. The fundamental spin-weights of $U_m$ are denoted by $\varpi_- = \varpi_m$ and $\varpi_+ = \varpi_{m+1}$, and similarly we write $\eta_- = \eta_{m-1}$ and $\eta_+ = \eta_m$ for the fundamental spin-weights of $M_m$. Let $(\sigma^\pm, V_{\sigma^\pm})$ be the two fundamental spin-representations of $M_m$ with highest weight $\eta_\pm$. There is only one fundamental spin-representation $(\tau, V_\tau)$ of $K_m$ with highest weight $\omega_m$. We denote by $(\pi_\varpi, V_\varpi)$ the finite-dimensional irreducible $\lag$-module (or $U_m$-module) of highest weight $\varpi$.

The embedding of $M_m \subset U_m$ is fixed such that $\eta_j$ corresponds to $\varpi_{j+1} - \varpi_1$ for $j < m-1$ and $\eta_\pm$ corresponds to $\varpi_\pm - \tfrac{1}{2}\varpi_1$. This identifies $\lat_\lam^* \subset \lat^*$, and the complement is given by $\lah^* = \bbC\varpi_1$. Dually, this fixes a one-dimensional subtorus $\lah \subset \lat$ for which $\lat = \lat_\lam\oplus\lah$ is an orthogonal decomposition with respect to the Killing form. The first fundamental weight $\varpi_1$ is now the only fundamental spherical weight; the only irreducible representations of $U_m$ that contain a $K_m$-fixed vector have as highest weight a multiple of $\varpi_1$. Similarly, the only irreducible representations of $U_m$ that contain the irreducible representation $\tau$ of $K_m$ have highest weight $\varpi_\pm + \ell\varpi_1$, with $\ell\in\bbZ_{\geq 0}$. Note that $\lah$ is the complexification of the Lie algebra of the Cartan torus $A_m \subset U_m$ that has the property that $U_m = K_m A_m K_m$.

Let $\rho_\lag$ and $\rho_\lam$ be the half sums of the positive roots for $(\lag, \lat)$ and $(\lam, \lat_\lam)$. The longest Weyl group elements in $W(\lag, \lat)$ and $W(\lam, \lat_\lam)$ are denoted by $w_{0, \lag}$ and $w_{0, \lam}$. The simple reflections in the Weyl group $W(\lag, \lat)$ that we will use are $s_{\eps_1\pm\eps_{m+1}}$.

The universal enveloping algebra of a Lie algebra $\las$ is denoted by $U(\las)$, and the center of $U(\las)$ is denoted by $\laZ_\las$. In the case $\las = \lag$, the action of the center $\laZ_\lag$ on irreducible $\lag$-modules is given by its central character. In particular, an element $Z \in \laZ_\lag$ acts on a Verma module $M(\lambda)$ by $\chi_\lambda(Z)$, where $\chi_\lambda$ denotes the central character of $M(\lambda)$. The central characters for Verma modules whose highest weights lie in the same Weyl orbit are equal, i.e., $\chi_\lambda = \chi_{w\lambda}$ for all $w \in W(\lag,\lat)$. The central character of $V_{\varpi_\pm + \ell\varpi_\sph}$ is
\begin{equation}\label{centralcharFDirrep}
\chi_{\varpi_\pm + \ell\varpi_1 + \rho_\lag}.
\end{equation}

Finally, we recall from \cite[(3.5.5)]{MR0929683} another important family of $\lag$-modules. The infinitesimal principal series representations $(H^{\sigma^\pm, \lambda\varpi_\sph})_{K_m}$ is the $(\lag, K_m)$-module that is obtained by taking the $K_m$-finite vectors of the principal series representations for the non-compact Cartan dual of compact symmetric pair $(U_m, K_m)$. Each $\lag$-module from the infinitesimal principal series admits a central character that we want to compare with \eqref{centralcharFDirrep}.

\begin{lemma}\label{lem:aflip}
Writing $w = s_{\eps_1 - \eps_{m+1}} s_{\eps_1 + \eps_{m+1}}$, we have
\begin{equation*}
w(\varpi_\pm + \mu\varpi_\sph + \rho_\lag) = \varpi_\mp - (\mu + 2m + 1)\varpi_\sph+\rho_\lag.
\end{equation*}
\end{lemma}

\begin{proof}
We have $w(\varpi_\sph) = -\varpi_\sph$, $w(\varpi_\pm) = \varpi_\mp - \varpi_\sph$ and $w(\varpi_j) = \varpi_j - 2\varpi_\sph$ for $2\leq j\leq m$. Hence $w(\rho_\lag) = \rho_\lag - 2m\varpi_\sph$ and $w(\varpi_\pm + \mu\varpi_\sph) = \varpi_\mp - (\mu + 1)\varpi_\sph$ from which the claim follows.
\end{proof}

\begin{lemma}\label{lem:centr_char}
The central character of $(H^{\sigma^\pm, \lambda\varpi_\sph})_{K_m}$ is given by $\chi_{\varpi_\pm + (\lambda - m - \tfrac{1}{2})\varpi_\sph + \rho_\lag}$.
\end{lemma}

\begin{proofempty}
The central character $\chi$ of $(H^{\sigma^\pm,\lambda\varpi_\sph})_{K_m}$ is $\chi_{w_{0,\lam}(\eta_\pm) + \lambda\varpi_\sph - \rho_\lam}$, see \cite[\S3.5.8,9]{MR0929683}. Note that $\rho_\lag = m \varpi_\sph + \rho_\lam$, which implies $w_{0,\lag}(\rho_\lam) = -\rho_\lag + m \varpi_\sph$. Combining this with $w_{0,\lag}(w_{0,\lam}(\eta_\pm)) = \eta_\mp$, we find
\begin{equation*}
\chi = \chi_{w_{0,\lag}({w_{0,\lam}(\eta_\pm) + \lambda\varpi_\sph - \rho_\lam})} = \chi_{\eta_\mp - (\lambda + m)\varpi_\sph + \rho_\lag}.
\end{equation*}
Using that $\eta_\pm = \varpi_\pm - \tfrac{1}{2}\varpi_\sph$ and with $w$ from Lemma \ref{lem:aflip}, we have
\pushQED{\qed}
\begin{equation*}
\chi = \chi_{w(\varpi_\mp - (\lambda + m + \tfrac{1}{2})\varpi_\sph + \rho_\lag)} = \chi_{\varpi_\pm + (\lambda - m - \tfrac{1}{2})\varpi_\sph + \rho_\lag}.\qedhere
\end{equation*}
\popQED
\end{proofempty}

\begin{corollary}\label{cor:subquotient_result}
The finite-dimensional irreducible $\lag$-module $V_{\varpi_\pm + \ell\varpi_\sph}$ is a subquotient of $(H^{\sigma^\pm, (\ell + m + \tfrac{1}{2})\varpi_\sph})_{K_m}$.
\end{corollary}

\begin{proof}
This follows from comparing the central characters from \eqref{centralcharFDirrep} and Lemma \ref{lem:centr_char}, see also \cite[Prop.3.2]{MR327978}.
\end{proof}

\begin{remark}
The principal series representations $H_{\sigma^\pm, \lambda}$, $\lambda \in \bbC$, that Camporesi and Pedon use in \cite{MR1814669} correspond to the $H^{\sigma^\pm, i\lambda\varpi_\sph}$ in our notation. The spherical functions $\Phi_\pm(\lambda)$ for the non-compact Cartan dual $(\Spin(2m+1,1), \Spin(2m+1))$ of $(U_m, K_m)$ have been investigated in \cite{MR1814669}. The reader should take note that the transformation behavior of the spherical functions in \cite{MR1814669} is inverted, and as to compare with our spherical functions the argument needs to be inverted. It follows from Corollary \ref{cor:subquotient_result} that $\Psi_\pm(\ell)(z)$ corresponds to $\Phi_\pm(\lambda(\ell))(z^{-1})$ with $\lambda(\ell) = -i(\ell + m + \tfrac{1}{2})$, which agrees with correspondence of the spectral parameters as established in \cite[Cor.3.5]{vHvP1}.
\end{remark}


\subsection{Lepowsky homomorphism}
We study $U_m$-modules $V$ that contain upon restriction to $K_m$ the irreducible representation $\tau$ of $K_m$. In this context, it is natural to consider the representations of the algebra $U(\lag)^\lak$, the centralizer of $\lak$ in $U(\lag)$, which acts naturally on $\Hom_\lak(V_\tau, V)$, the $\tau$-isotypical component of $V$. If $V$ is an irreducible $U_m$-module that contains $V_\tau$ upon restriction to $K_m$, then $U(\lag)^\lak$ acts by a character because the multiplicity is at most one. This character can be understood as follows. 

The Lepowsky homomorphism
\begin{equation*}
L : U(\lag)^\lak \to S(\lah) \otimes U(\lak)^\lam
\end{equation*}
is a projection onto a direct summand of $U(\lag)$, see \cite[\S3.5]{MR0929683}. This projection is in spirit similar to the projection $U(\lag) \to S(\lat)$ in the proof of Harish-Chandra's isomorphism $\laZ_\lag\to S(\lat)^W$, see \cite[Thm.3.2.3]{MR0929683}. The action of $u\in U(\lag)^\lak$ on $\Hom_\lak(V_\tau, (H^{\sigma^\pm, \lambda\varpi_\sph})_{K_m})$ is given by the action of $L(u)((\lambda + m)\varpi_\sph)$ on $\Hom_\lam(V_{\sigma^\pm}, V_\tau)$, see \cite[(3.5.7,8)]{MR0929683}. To indicate that we are interested how $u \in  U(\lag)^\lak$ acts on $\tau$-isotypical components of representations with label $\sigma^\pm$, we will write $L_{\tau,\sigma^\pm}(u)$, which can be viewed as an element of $S(\lah)$ because $L_{\tau,\sigma^\pm}(u)(\mu\varpi_\sph)$ acts by a scalar. For $u \in U(\lag)^\lak$, we define $L_\tau(u) \in S(\lah)\otimes\bbM_2$ by
\begin{equation*}
L_\tau(u)(\lambda\varpi_\sph) =
\begin{pmatrix}
L_{\tau, \sigma^-}(u)(\lambda\varpi_\sph) & 0 \\
0 & L_{\tau, \sigma^+}(u)(\lambda\varpi_\sph)
\end{pmatrix}.
\end{equation*}
To get control over these scalars we use that the algebra $U(\lag)^\lak$ is in this case isomorphic to $\laZ_\lag \otimes \laZ_\lak$ in our particular example, see \cite[Satz 2.3]{MR1048234}. 

\begin{proposition}\label{prop:invariance}
For $u\in U(\lag)^\lak$, we have
\begin{equation*}
\Rad_\tau(u)\Psi(\ell) = \Psi(\ell)L_\tau(u)((\ell+2m+\tfrac{1}{2})\varpi_\sph).
\end{equation*}
Moreover, the function
\begin{equation*}
\ell + \tfrac{1}{2}\mapsto L_\tau(u)((\ell+2m+\tfrac{1}{2})\varpi_\sph)
\end{equation*}
is invariant for the action \eqref{eqn:s_*_action} of $\bbZ_2$.
\end{proposition}
\begin{proof}
The eigenvalue property of $L_\tau(u)$ follows from its definition. For the $\bbZ_2$-invariance, note that for $u \in \laZ_\lag$ we have $L_{\tau, \sigma^\pm}(u)((\ell + 2m + \tfrac{1}{2})\varpi_\sph) = \chi_{\varpi_\pm + \ell\varpi_\sph + \rho_\lag}(u)$. The application of $s_*$ changes $\ell + \tfrac{1}{2} \mapsto -\ell - 2m - \tfrac{1}{2}$ and it becomes $\chi_{\varpi_\pm - (\ell + 2m + 1)\varpi_\sph + \rho_\lag}(u)$, which is equal to $\chi_{\varpi_\mp +\ell\varpi_\sph + \rho_\lag}(u)$ by Lemma \ref{lem:aflip} and in turn equal to $L_{\tau, \sigma^\mp}(u)((\ell + 2m + \tfrac{1}{2})\varpi_\sph)$.
\end{proof}

We define another action of $\bbZ_2$ on $\lah^*$, given by $s^*(\lambda\varpi_\sph) = s(\lambda\varpi_\sph - m\varpi_\sph) + m\varpi_\sph$, and we denote by $S(\lah)\otimes\bbC^2$ the algebra of $\bbM_2$-valued polynomials with values in the diagonal matrices. Let $\left(S(\lah) \otimes \bbC^2\right)^{\bbZ_2^*}$ denote the space of invariants for the action
\begin{equation*}
s(p\otimes A)(\lambda\varpi_\sph) = p((s^*)^{-1}\lambda\varpi_\sph) \otimes \left(\pi(s)A\pi(s)^{-1}\right).    
\end{equation*}
As a consequence of Proposition \ref{prop:invariance}, $L_\tau$ maps into $\left(S(\lah) \otimes \bbC^2\right)^{\bbZ_2^*}$. Using invariant theory, one can argue that $L_\tau$ is also surjective, however, we give a direct argument using the radial part map.

Let $\frakR'$ be the unital algebra of $\bbC$ generated by the functions $(1-z^{-2})^{-1}$ and $z^{\pm1}$. Then $\frakR'$ is a subalgebra of the algebra of meromorphic functions on $A_m$. There exists a homomorphism
\begin{equation*}
\Rad_\tau : U(\lag)^\lak \to \frakR' \otimes S(\lah) \otimes \bbM_2
\end{equation*}
such that $\Rad_\tau(u)$ is a differential operator for which the spherical functions, restricted to $A_m$, are eigenfunctions, see \cite[Thm.3.1]{MR683007}. This implies that $\Rad_\tau(u)$ is $\bbZ_2$-invariant. Let $I_\tau \subset U(\lak)$ be the kernel of the representation $U(\lak) \to \End(V_\tau)$. Then both $L_\tau$ and $\Rad_\tau$ factor through the commutative algebra
\begin{equation*}
\bbD(\tau) = U(\lag)^\lak / \left(U(\lag)^\lak \cap U(\lag)I_\tau\right).
\end{equation*}

Let $\calD(\tau)=\bbC[R_m]$ denote the image of $\Rad_\tau$. Then the diagram
\begin{equation*}
\xymatrix{\bbD(\tau)\ar[rr]^{\overline{\Rad}_\tau}\ar[rd]_{\overline{L}_\tau}&&\calD(\tau)
\ar[ld]\\ & \left( S(\lah) \otimes \bbC^2 \right)^{\bbZ_2^*}}
\end{equation*}
implicitly defines an isomorphism $\calD(\tau) \to \left(S(\lah) \otimes \bbC^2 \right)^{\bbZ_2^*}$ that can be understood as an eigenvalue map; evaluate in $(\ell + 2m + \tfrac{1}{2})\varpi_\sph$ for the eigenvalue of $\Psi(\ell)$. Since $R_m$ has eigenvalues that are linear in $\ell$, we see that $L_\tau$ is surjective. 

We conclude that eigenvalues of the full spherical functions $\Psi(\ell)$ are invariant $\bbM_2$-valued polynomials in $\ell+\tfrac12$, with values taken in the diagonal matrices in $\bbM_2$. The components of this polynomial are given by the Lepowsky homomorphism $L$, a projection map of $U(\lag)$ onto a direct summand. 

From this eigenvalue map we can also obtain the polynomial that we need to evaluate $R_m$ in to get back the operator in $\calD(\tau)$ that we started with. In this sense we recover the Harish-Chandra map from Section \ref{S:Eigenvalue_map} in the context of representation theory. This result is also analogous to the interpretation of the Harish-Chandra homomorphism $U(\lag)^\lak \mapsto S(\lah)$ \cite[Thm.5.18]{MR1790156} that encodes the eigenvalues for the invariant differential operators acting on the zonal spherical functions. 


\appendix
\section{}\label{app:a}

We abbreviate the action of $s$ on $f \in \bbC[z^{\pm 1}]$ by $f^s = s(f)$.

\begin{lemma}\label{lem:derivative_of_f1_f2}
Let $f \in \bbC[z^{\pm 1}]$ and recall from Lemma \ref{lem:steinberg_basis} that $f = f_1 + f_2 z$ with
\begin{equation*}
f_1 = \frac{z f^s - z^{-1} f}{z - z^{-1}}, \quad f_2 = \frac{f - f^s}{z - z^{-1}}.
\end{equation*}
Then
\begin{align*}
z\partial_z(f_1) &= \frac{z f^s + z^2\partial_z(f^s) + z^{-1} f - \partial_z(f)}{z - z^{-1}} - \frac{(zf^s - z^{-1} f)(z + z^{-1})}{(z - z^{-1})^2}, \\
z\partial_z(f_2) &= \frac{z\partial_z(f) - z\partial_z(f^s)}{z - z^{-1}} - \frac{(f - f^s)(z + z^{-1})}{(z - z^{-1})^2}.
\end{align*}
\end{lemma}

\begin{proof}
This follows from a small computation.
\end{proof}

\begin{proofempty}[ of Proposition \ref{prop:nonsym_shifts_expl}]
Let $f \in \bbC[z^{\pm 1}]$. Note that $\widehat{G}_+^U(k) = G_+(k)I_2$, and using Lemma \ref{lem:derivative_of_f1_f2} we compute
\begin{equation*}
\calG_+(k)(f) = \frac{1}{z - z^{-1}}(z\partial_z(f_1) + (z\partial_z(f_2))z) = G_+(k)(f) - \frac{z (f - f^s)}{(z - z^{-1})^2}.
\end{equation*}
Similarly, note that
\begin{equation*}
\widehat{G}^U_-(k) = \left( G_-(k) + (z + z^{-1}) \right)I_2 - 2
\begin{pmatrix}
0 & 1 \\ 1 & 0
\end{pmatrix}, \quad \begin{pmatrix} 0 & 1 \\ 1 & 0 \end{pmatrix} \underline{f} = \underline{zf^s},
\end{equation*}
which imply that
\begin{align*}
\calG_-(k)(f) &= (z - z^{-1})(z\partial_z(f_1) + (z\partial_z(f_2))z) + ((k_1 + 2k_2)(z + z^{-1}) + 2k_1)f - 2zf^s \\
&= G_-(k)(f) + z^{-1} f - z f^s.
\end{align*}
For the last two cases, we have that
\begin{equation*}
\widehat{E}^U_\pm(k) = E_\pm(k)I_2 + \begin{pmatrix} 0 & \pm1 \\ 0 & 1 \end{pmatrix}, \quad \begin{pmatrix} 0 & \pm1 \\ 0 & 1 \end{pmatrix} \underline{f} = \underline{(z \pm 1)\frac{f - f^s}{z - z^{-1}}},
\end{equation*}
which give that
\begin{align*}
\calE_+(k)(f) &= E_+(k)(f) -\frac{1 + z^{-1}}{1 - z^{-1}} \frac{z (f - f^s)}{z - z^{-1}} + (z + 1) \frac{f - f^s}{z - z^{-1}} \\
&= E_+(k)(f) -\frac{1 + z^{-1}}{1 - z^{-1}} \frac{f - f^s}{z - z^{-1}},
\end{align*}
and
\pushQED{\qed}
\begin{align*}
\calE_-(k)(f) &= E_-(k)(f) - \frac{1 - z^{-1}}{1 + z^{-1}} \frac{z (f - f^s)}{z - z^{-1}} + (z - 1) \frac{f - f^s}{z - z^{-1}} \\
&= E_-(k)(f) + \frac{1 - z^{-1}}{1 + z^{-1}} \frac{f - f^s}{z - z^{-1}}. \qedhere
\end{align*}
\popQED
\end{proofempty}

\bibliography{MvH-and-MvP-Non-sym-BC1-part-2}{}
\bibliographystyle{plain}

\end{document}